\documentclass[11pt]{article}
\usepackage{graphicx} 
\usepackage[T1]{fontenc} 
\usepackage[utf8]{inputenc}
\usepackage{amsmath,amssymb}
\usepackage[table,dvipsnames]{xcolor}
\usepackage{amsthm}
\usepackage{graphicx}
\usepackage{tikz}
\usepackage{dsfont}
\usepackage{nicefrac}
\usepackage{booktabs} 
\usepackage{pgfplots}
\usepackage{enumitem}
\usepackage{url}
\usepackage{algorithm}
\usepackage{algorithmic} 
\usepackage[normalem]{ulem}

\usetikzlibrary{shapes, positioning, quotes}
\pgfplotsset{compat=1.18} 

\DeclareMathOperator{\Div}{div}
\newcommand{\mbf}[1]{\boldsymbol{#1}}

\usepackage{multirow}

\newcommand{\bb}[1]{\mathbb{#1}}
\newcommand{\vertiii}[1]{{\left\vert\kern-0.25ex\left\vert\kern-0.25ex\left\vert #1 
    \right\vert\kern-0.25ex\right\vert\kern-0.25ex\right\vert}}

\usetikzlibrary{decorations.pathreplacing,calc}
\newcommand{\tikzmark}[1]{\tikz[overlay,remember picture] \node (#1) {};}

\newcommand*{\AddNote}[5]{%
    \begin{tikzpicture}[overlay, remember picture]
        \draw [decoration={brace,amplitude=0.8em},decorate,very thick,#4]
            ($(#3)!(#1.north)!($(#3)-(0,1)$)$) --  
            ($(#3)!(#2.south)!($(#3)-(0,1)$)$)
                node [align=left, text width=8cm, pos=0.5, anchor=west] {#5};
    \end{tikzpicture}
}%

\newlength\myindent
\usepackage{bm}

\newcommand\y{\cellcolor{red!20}}

\usepackage{a4wide}
\theoremstyle{remark}
\newtheorem*{remark}{Remark}

\title{$\varphi$-FEM-FNO: a new approach to train a Neural Operator as a fast PDE solver for variable geometries}

\author{Michel Duprez\footnote{
MIMESIS team, Inria de l'Université de Lorraine, MLMS team, Universit\'e de Strasbourg, 2 Rue Marie Hamm, 67000 Strasbourg, France, \texttt{michel.duprez@inria.fr}} ,
Vanessa Lleras\footnote{IMAG, Univ Montpellier, CNRS UMR 5149, 499-554 Rue du Truel, 34090 Montpellier, France,
\texttt{vanessa.lleras@umontpellier.fr}} ,
Alexei Lozinski\footnote{Université de Franche-Comté, Laboratoire de mathématiques de Besançon, UMR~CNRS~6623, 16 route de Gray, 25030 Besançon Cedex, France, 
\texttt{alexei.lozinski@univ-fcomte.fr}} , \\ 
Vincent Vigon\footnote{Institut de Recherche Mathématique Avancée, UMR 7501, Université de Strasbourg et CNRS, Tonus team, Inria de l'Université de Lorraine, 7 rue René Descartes, 67000 Strasbourg, France, 
\texttt{vincent.vigon@math.unistra.fr}} 
\ and 
Killian Vuillemot\footnote{IMAG, Univ Montpellier, CNRS UMR 5149, 499-554 Rue du Truel, 34090 Montpellier, France. MIMESIS team, Inria de l'Université de Lorraine, MLMS team, Universit\'e de Strasbourg, 2 Rue Marie Hamm, 67000 Strasbourg, France, 
\texttt{killian.vuillemot@umontpellier.fr}}
}

\begin{document}






\maketitle

\begin{abstract}
    In this paper, we propose a way to solve partial differential equations (PDEs) by combining machine learning techniques and the finite element method called $\varphi$-FEM. For that, we use the Fourier Neural Operator (FNO), a learning mapping operator. The purpose of this paper is to provide numerical evidence to show the effectiveness of this technique. We will focus here on the resolution of two equations: the Poisson-Dirichlet equation and the non-linear elasticity equations.
    The key idea of our method is to address the challenging scenario of varying domains, where each problem is solved on a different geometry. The considered domains are defined by level-set functions due to the use of the $\varphi$-FEM approach. We will first recall the idea of $\varphi$-FEM and of the Fourier Neural Operator. Then, we will explain how to combine these two methods. We will finally illustrate the efficiency of this combination with some numerical results on three test cases. In addition, in the last test case, we propose a new numerical scheme for hyperelastic materials following the $\varphi$-FEM paradigm.
\end{abstract}



\section{Introduction}

Finite Element Method (FEM) is one of the most popular approaches to approximate the solutions of Partial Differential Equations (PDE) arising in engineering, physics, biology, and other applications (see e.g. \cite{ern2004theory}).
It is important to solve them quickly (sometimes in real-time) with good accuracy. There have been numerous attempts to achieve this using machine learning-based (ML-based) methods. They can be split into two groups :
\begin{enumerate}
    \item \textbf{Physics-inspired approaches: }
          ML-based methods can be used as an approximation ansatz and approximate the solution of PDEs by minimizing the residual or the associated energy of the PDEs and the distance to some observations, without ever using traditional approximation by FEM or similar. The most popular member of this class of methods is PINNs \cite{PINNs}, but one can also cite  Deep Galerkin \cite{sirignano2018dgm} and Deep Ritz methods \cite{yu2018deep}. Despite the initial promise, there is now abundant numerical evidence that these methods do not outperform the classical FEM in terms of solution time and accuracy, see for instance a recent study in \cite{grossmann2023can}. It seems that these methods cannot thus be considered as good candidates for real-time realistic computations.
    \item \textbf{Classical solver as database: }
          Classical FEM (or similar) is used to obtain a "database" of solutions for a collection of representative parameter values that are used to train a neural network to learn the mapping linking the parameters to the solution. This step is computationally expensive and is done in the preparatory stage (offline). The expected outcome is that one can use the trained network to obtain the solution for any given parameters almost instantaneously (online). Examples include U-Net (see e.g. \cite{ronneberger2015u}), Graph Neural Operator \cite{li2020neural}, DeepOnet \cite{lu2019deeponet} and Fourier Neural Operator (FNO) \cite{paper_FNO,FNO_general_geometries}.
\end{enumerate}
In our article, we focus on FNO as the method that showed a superior cost-accuracy tradeoff over the others (see \cite{paper_FNO}).
The issue with FNO is that it needs Cartesian grids to perform discrete fast Fourier transform, and the initial implementation was thus limited to problems posed on rectangular boxes.  There have been attempts to adapt FNO to general geometries, cf. Geo-FNO \cite{FNO_general_geometries} where the irregular input domain is deformed into a uniform latent mesh on which the FFT can be applied. In our article, we propose an alternative approach: we treat the geometry, given by the level-set function, as one of the inputs of the network alongside the other data of the problem while using a Cartesian grid without deforming it. Incidentally, this viewpoint of treating the geometry (i.e. thanks to the level-set function) together with the data to construct an approximation using a simple (ex. Cartesian) grid was also the starting point to develop $\varphi$-FEM. It is thus natural to combine $\varphi$-FEM (at the offline training stage) with FNO. As a bonus, this combination, which we will call $\varphi$-FEM-FNO, allows us to avoid the interpolation errors from a body-fitted mesh to a Cartesian one, which would be inevitable if we used a traditional  FEM for training.

This paper aims to illustrate the efficiency of our approach $\varphi$-FEM-FNO, in the case of complex and varying domains
for the Poisson-Dirichlet problem:
\begin{equation}\label{eq:governing_poisson}
    \begin{cases}
        - \Delta u & = f \,, \quad \text{ in } \Omega\,, \\
        \hfill u   & = g \,, \quad \text{ on } \Gamma\,,
    \end{cases}
\end{equation}
and for the non-linear elasticity equations:
\begin{equation}
    \label{eq:governing_hyperelast}
    \begin{cases}
        - \Div \mbf{P}(\mbf u)              & = \mbf{f}\,, \quad \hfill \text{in } \Omega \,,        \\
        \hfill \mbf u                       & = \mbf u_D\,, \quad \hfill \text{on } \Gamma_D\,,      \\
        \hfill \mbf{P}(\mbf u) \cdot \mbf n & = \mbf{t}\,, \quad \hfill \quad \text{on } \Gamma_N\,,
    \end{cases}
\end{equation}
where $\Omega$ is a connected domain of $\bb{R}^d,d=1,2,3$ and $\Gamma$ its boundary, with $\Gamma = \Gamma_D \cup \Gamma_N$ and $\Gamma_D \cap \Gamma_N = \emptyset$ in \eqref{eq:governing_hyperelast}.

Our contributions are the following:
\begin{itemize}
    \item We propose a new machine learning approach called $\varphi$-FEM-FNO which takes as input the parameters of the PDE and the geometry of the domain encoded by a level-set function $\varphi$, and gives as output an approximation of the solution of the PDE. In comparison to \cite{FNO_general_geometries}, we do not need a transformation between the geometry and the unit square. Our approach is thus simpler and results in a lighter and smaller operator to train.
    \item In addition, we highlight in Fig.~\ref{fig:compare_methods_test_case_1} that our approach has a better accuracy/CPU-time ratio than Geo-FNO, which has been compared with other techniques in \cite{FNO_general_geometries}, and than an FNO trained using standard FEM solutions interpolated on cartesian grids {or a UNet trained using the same approach as for $\varphi$-FEM-FNO}.
    \item We introduce and validate a new $\varphi$-FEM scheme to solve hyperelastic problems on complex
          geometries and propose a combination of this scheme with the $\varphi$-FEM-FNO approach.
    \item We also introduce an ML-based method constructed using a classical finite element solver and a level-set description of the domains. This method we called Standard-FEM-FNO is very simple but also provides quite interesting results compared to our main method $\varphi$-FEM-FNO.
\end{itemize}

The paper will be divided into three parts. In Section \ref{sec:descr}, we will first describe the two methods used: $\varphi$-FEM and FNO. We will then present in Section \ref{sec:main idea} our idea to combine them. For readability, these two sections will be devoted to the treatment of the Poisson-Dirichlet equation \eqref{eq:governing_poisson}.
Finally, we will illustrate the efficiency of the method with numerical results in Section \ref{sec:num} and give some conclusions in Section \ref{sec:concl}. In the end, we detail in the first appendix the standardization operator, and we present the optimizer algorithm in the second appendix.

\section{Description of the methods}\label{sec:descr}
In the rest of the manuscript, $\Omega$ is a domain of dimension $2$ included in $[0,1]^2$.
Moreover, we focus in this section on the resolution of the equation \eqref{eq:governing_poisson}.
\subsection{Overview}

Our idea is to build a neural network that will be an approximation of the operator mapping the data $f$, $g$, and the geometry to the solution of \eqref{eq:governing_poisson}. We want the output to be obtained with good accuracy and a low computational time. The objective is to train this neural network using synthetic data generated by a discrete solver of PDE. The neural network and the discrete solver must be chosen to perform independently of each other, and must also be compatible.

As a discrete solver, we choose $\varphi$-FEM \cite{phifem} which is a finite element method with an immersed boundary approach using a level-set function to describe the geometry of the domain.
The optimal convergence of $\varphi$-FEM  has been previously proven theoretically and numerically for the Poisson equation with Dirichlet boundary conditions \cite{phifem}, with Neumann boundary conditions \cite{phiFEM2}, for the Stokes problem  \cite{phifemstokes} and for the Heat-Dirichlet equation \cite{phiFEM_heat}.
Moreover, in \cite{cotin:hal-03372733}, the efficiency of the method {compared to the continuous Lagrange FEM approach on conformal meshes}  has been illustrated numerically on multiple examples in the case of linear elasticity.

As neural network, we have decided to use the Fourier Neural Operator (FNO), introduced in \cite{paper_FNO} and \cite{Neural_operator}.
The FNO relies on an iterative architecture proposed in \cite{li2020neural}. An advantage of FNO is that it takes a step size much bigger than is allowed in numerical methods.
In the case of the approximation of the PDE solution, the authors of \cite{paper_FNO} have illustrated that FNO has better performance than the classical Reduced Basis Method (using a POD basis) \cite{devore2017theoretical},
a Fully Convolution Networks \cite{zhu2018bayesian}, an operator method using PCA as an autoencoder
on both the input and output data and interpolating the latent spaces with a neural network \cite{bhattacharya2021model}, the original graph neural operator \cite{li2020neural}, the
multipole graph neural operator \cite{li2020multipole}, a neural operator method based on the
low-rank decomposition of a kernel similar to the unstacked DeepONet proposed in \cite{lu2019deeponet}. They have also illustrated that FNO outperforms
a ResNet (18 layers of 2-d convolution with residual connections) \cite{he2016deep},
a U-Net \cite{ronneberger2015u} and  TF-Net \cite{wang2020towards}.
Furthermore, the training {of FNOs} can be done on many PDEs with the same underlying architecture.

Moreover, these two methods are compatible since $\varphi$-FEM is a precise non-conforming finite element method, that can be used on cartesian grids, as required by the FNO that will be used.

In the next subsections, we will describe $\varphi$-FEM and FNO to solve the equation \eqref{eq:governing_poisson}. We introduce in Table~\ref{table:notations} notations that will be used in the rest of the manuscript.
\begin{table}
    \begin{center}
        \begin{tabular}{ | c | c | c |}
            \hline
             & \textbf{Notation}                 & \textbf{Meaning}                                                          \\ \hline
            \multirow{6}*{\rotatebox{90}{FNO}}
             & $\theta$                          & Set of trainable parameters                                               \\
             & $\mathcal{G}_{\theta}$            & Operator mapping the input to the solutions                               \\
             & $\mathcal{G}^{\dagger}$           & Ground truth operator mapping the solutions                               \\
             & $\mathcal{F}$, $\mathcal{F}^{-1}$ & Discrete Fourier and inverse Fourier transform                            \\
             & $W^{\mathcal{C}^l_\theta}$        & Linear transformation applied on lower Fourier modes                      \\
             & $\mathcal{C}^l_\theta$            & Convolution layer                                                         \\
             & $\mathcal{B}^l_\theta$            & Linear transformations applied on the spatial domain                      \\
             & $P_\theta$, $Q_\theta$            & Transformations between high dimension channel space and original space   \\
             & $N$, $N^{-1}$                     & Standardization and unstandardization operators                           \\
             & $\sigma$                          & Non linear activation function                                            \\ \hline
            \multirow{6}*{\rotatebox{90}{$\varphi$-FEM}}
             & $\varphi$                         & Level-set function defining the domain $\Omega$ and its boundary $\Gamma$ \\
             & $\mathcal{O}$                     & Box $[0,1]^2$                                                             \\
             & $\mathcal{T}_h$                   & $\varphi$-FEM computational mesh                                          \\
             & $\mathcal{T}_h^\Gamma$            & Set of cells of $\mathcal{T}_h$ cut by the boundary                       \\
             & $\mathcal{F}_h^\Gamma$            & Set of internal facets of $\mathcal{T}_h^\Gamma$                          \\
             & $\sigma_D$                        & Stabilisation parameter                                                   \\
            \hline
        \end{tabular}
        \\ \caption{Notations table.}\label{table:notations}
    \end{center}
\end{table}

\subsection{Description of $\varphi$-FEM}

Let us first briefly describe the $\varphi$-FEM method introduced in \cite{phifem} to solve \eqref{eq:governing_poisson}. We will skip many theoretical aspects but refer the reader to \cite{phifem} for more details.
We suppose that the domain $\Omega$ is included in the box $\mathcal{O}=[0,1]^2\subset\mathbb{R}^2$ and is given by a level-set function $\varphi$ such that:

\begin{equation*}
    \Omega:= \{ \varphi < 0 \} \qquad \text{ and } \qquad \Gamma:= \{ \varphi = 0 \}\,,
\end{equation*}
where $\Gamma$ is the boundary of $\Omega$.

\begin{figure}
    \centering
    \includegraphics[width=0.42\textwidth]{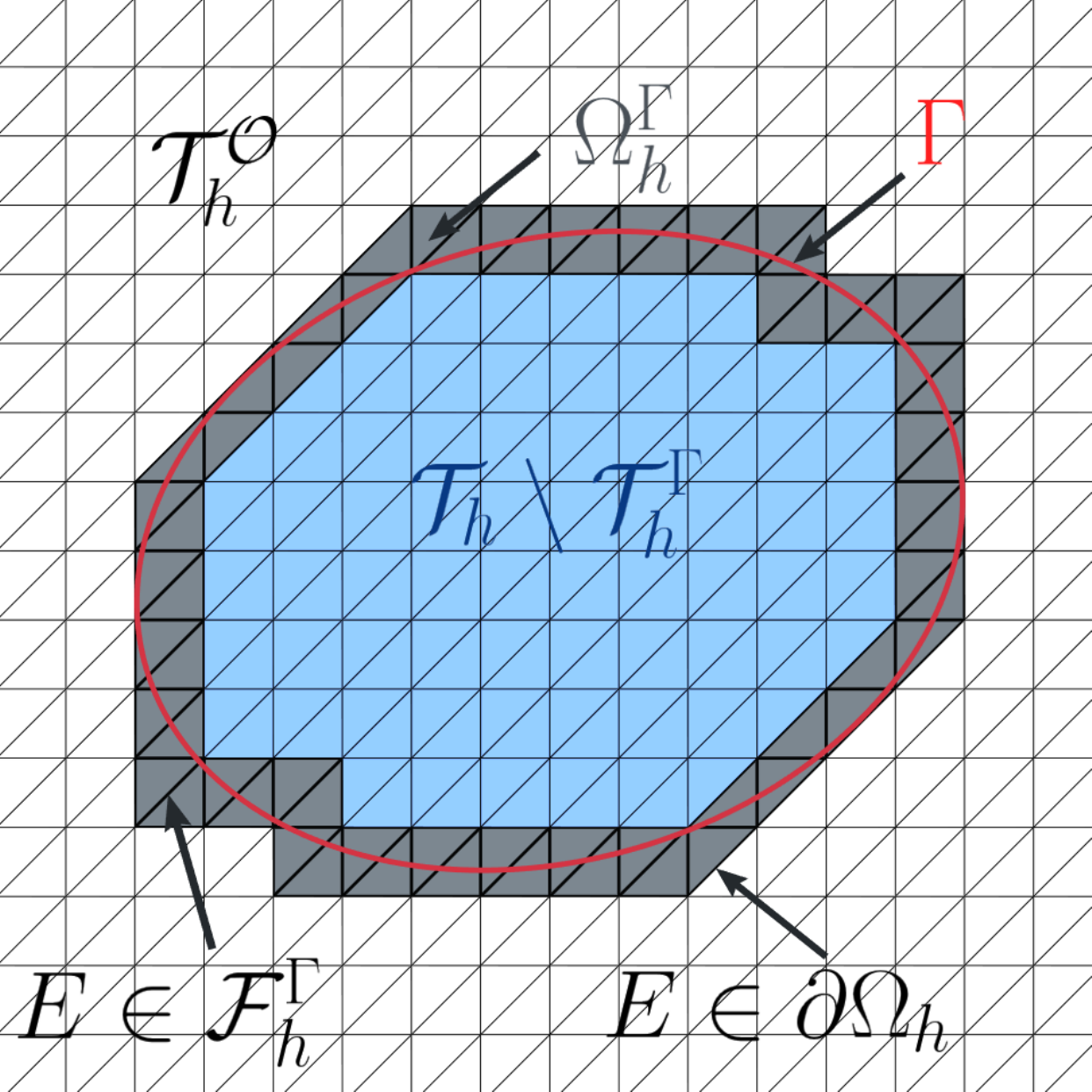}
    \quad
    \includegraphics[width=0.42\textwidth]{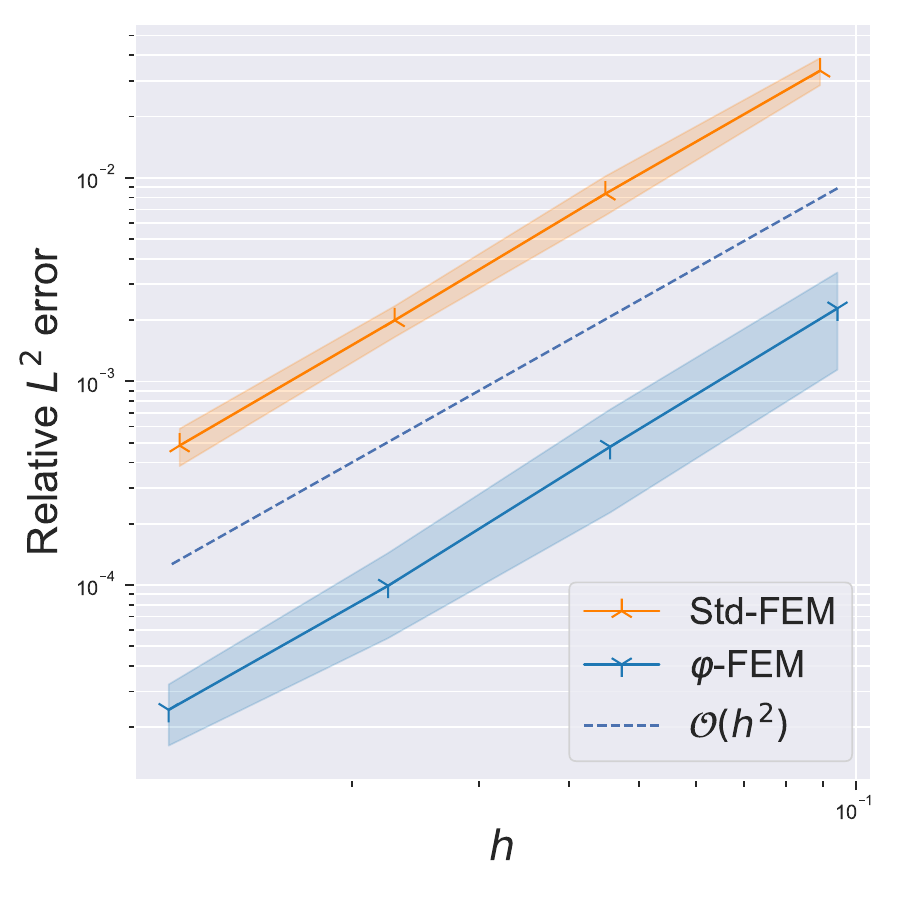}
    \caption{
    Left: example of $\varphi$-FEM meshes. In red, the exact boundary $\Gamma$ of an ellipse $\Omega$, in white $\mathcal{T}_h^{\mathcal{O}}$, in gray $\mathcal{T}_h^{\Gamma}$ and in blue, $\mathcal{T}_h \setminus \mathcal{T}_h^{\Gamma}$. Right: Convergence curves of $\varphi$-FEM and a standard finite element method, to solve \eqref{eq:governing_poisson} for 5 combinations of domain, force, and boundary conditions.    }
    \label{fig:domain_submeshes}
\end{figure}

Let $\mathcal{T}_h^{\mathcal{O}}$ be a triangular cartesian mesh of $\mathcal{O} =[0,1]^2$ composed of $n_x -1$ and $n_y-1$ squares divided into triangular cells in its width and its height, of sizes $h$.
Denoting by $\varphi_h$ the Lagrange interpolation of $\varphi$ on $\mathcal{T}_h^{\mathcal{O}}$, we consider the submesh $\mathcal{T}_h$ of $\mathcal{T}_h^{\mathcal{O}}$, called the computational mesh, composed of the cells of $\mathcal{T}_h^{\mathcal{O}}$ intersecting the domain $\{ \varphi_h < 0\}$, i.e.

\[
    \mathcal{T}_h:= \left\{T \in \mathcal{T}_h^{\mathcal{O}} : \ T \cap \{ \varphi_h < 0 \} \neq \emptyset  \right\}\,.
\]
We also introduce the submesh $\mathcal{T}_h^\Gamma$ containing the cells cut by the approximate boundary ($\{ \varphi_h = 0 \}$), i.e.
\[
    \mathcal{T}_h^\Gamma:= \left\{T \in \mathcal{T}_h \: \ T \cap \{ \varphi_h = 0 \} \neq \emptyset  \right\}\,.
\]
We denote by $\Omega_h$ and $\Omega_h^\Gamma$ the domains occupied by $\mathcal{T}_h$ and $\mathcal{T}_h^\Gamma$, respectively, and by $\partial\Omega_h$ the boundary of $\Omega_h$ (different from $\Gamma_h = \{ \varphi_h = 0 \}$).
All the meshes are illustrated for a specific domain in Fig.~\ref{fig:domain_submeshes} (left).
Finally, we need to introduce a set of facets containing all the internal faces of the mesh $\mathcal{T}_h^\Gamma$, i.e., the faces of $\mathcal{T}_h^\Gamma \setminus \partial \mathcal{T}_h$. Referring to Fig.~\ref{fig:domain_submeshes} (left), these faces are the ones of the gray cells except the ones common to a gray cell and a white cell. We will denote by $\mathcal{F}_h^\Gamma$ this set, defined by
\[
    \mathcal{F}_h^\Gamma:= \{ F \text{ (an internal facet of } \mathcal{T}_h \text{) such that } \exists \ T \in \mathcal{T}_h: T \cap \Gamma_h \neq \emptyset \text{ and } F \in \partial T \} \,.
\]

Let $k\geqslant 1$ be an integer. We define the finite element space
\[
    V_h^{(k)}:= \{ v_h \in H^1(\Omega_h) \: \ v_h|_T \in \bb{P}_k(T) \ \forall \ T \in \mathcal{T}_h \} \,.
\]

We now introduce the considered $\varphi$-FEM formulation of system \eqref{eq:governing_poisson} (see \cite{phifem}):
Find $w_h \in V_h^{(k)}$ such that, for all $s_h \in V_h^{(k)}$, denoting $u_h=\varphi_hw_h+g_h$ and $v_h=\varphi_hs_h$,
\[
    \int_{\Omega_h} \nabla u_h \cdot \nabla v_h  -\int_{\partial\Omega_h} \frac{\partial u_h}{\partial n}  v_h + G_h(u_h,v_h)= \int_{\Omega_h} f_h v_h + G_h^{rhs}(v_h) \,,
\]
where $g_h$, $f_h$ are some Lagrange interpolations of $g$ and $f$, respectively,
\begin{align*}
    G_h(u,v) & = \sigma_D h \sum_{E \in \mathcal{F}_h^\Gamma} \int_E \left[ \partial_n u \right]\left[ \partial_n v \right]  + \sigma_D h^2 \sum_{T \in \mathcal{T}_h^\Gamma}\int_T \Delta u \Delta v \,,
\end{align*}
and
\begin{align*}
    G_h^{rhs}(v) & = - \sigma_D h^2 \sum_{T \in \mathcal{T}_h^\Gamma}\int_T f_h \Delta v \,.
\end{align*}
The brackets in $G_h$ stand for the jump over the facets of $\mathcal{F}_h^\Gamma$ , $\partial_n u$ stands for the normal derivative of $u$ and $\sigma_D > 0 $ is a $h$-independent parameter.
Figure~\ref{fig:domain_submeshes} (right) illustrates the convergence of $\varphi$-FEM (measured in relative $L^2$ norm) and a standard finite element method in solving \eqref{eq:governing_poisson} for five combinations of domain $\Omega$, force $f$, and boundary condition $g$. Both methods exhibit the same order of convergence; however, the error of $\varphi$-FEM is significantly smaller than the one of the classical method. This observation strongly supports the choice of $\varphi$-FEM as the finite element solver.

\subsection{The ``ground truth'' operator}

In the rest of the manuscript, without specific mention, $f_h$, $g_h$, $\varphi_h$, $u_h$ and $w_h$ will represent the matrices of $\mathbb{R}^{n_x\times n_y}$ associated to these $\mathbb{P}_1$-functions composed for each index $i=0,\ldots,n_x-1$, $j=0,\ldots,n_y-1$, of the values of the evaluation or an extrapolation in $V_h^{\mathcal{O}}$ of these functions at the node of the mesh $\mathcal{T}_h^{\mathcal{O}}$ of coordinate $(x_i,y_j)$, with $x_i:=i/(n_x-1)$, $y_j:=j/(n_y-1)$, where
\[
    V_h^{\mathcal{O}}:= \{ v_h \in H^1(\mathcal{O}) \: \ v_h|_T \in \bb{P}_k(T) \ \forall \ T \in \mathcal{T}_h^{\mathcal{O}} \} \,.
\]

In the tradition of FNO literature, the FNO will approximate an operator called the ``ground truth operator'' which is denoted by $\mathcal{G}^{\dagger}$. In our case, $\mathcal{G}^{\dagger}$ will be the operator mapping the data $f_h$, $g_h$, and the geometry given by the level-set $\varphi_h$ to the $\varphi$-FEM approximated solution $w_h$. More precisely, $\mathcal{G}^{\dagger}$ will be defined as follows:
\begin{equation}\label{GdaggerVh}
    \begin{array}{rccl}
        \mathcal{G}^{\dagger}: & \mathbb R^{n_x\times n_y \times 3} & \to     & \mathbb R^{n_x\times n_y\times 1} \\
                               & (f_h,\varphi_h,g_h)                & \mapsto & w_h\,,
    \end{array}
\end{equation}
where $w_h$ is the $\varphi$-FEM solution associated to $f_h$, $g_h$ and $\varphi_h$.
Note that $w_h$ is simply extrapolated by $0$ outside $\Omega_h$, with no impact on the FNO since these $0$ values are not seen in the loss defined below. {Hence, we do not necessarily need explicit expressions of $f$ and $g$.} In practice, this extrapolation will be done by DOLFINx (\cite{DOLFINx,basix1,basix2,uflX}).

\subsection{Architecture of the FNO}

We will now introduce a few essential points to understand the architecture of the FNO.
We refer the reader to \cite{paper_FNO,Neural_operator,FNO_general_geometries} for detailed explanations on the FNO or \cite{li2020neural} for more details about neural operators.

\begin{figure}
    \centering
    \includegraphics[width=0.9\textwidth]{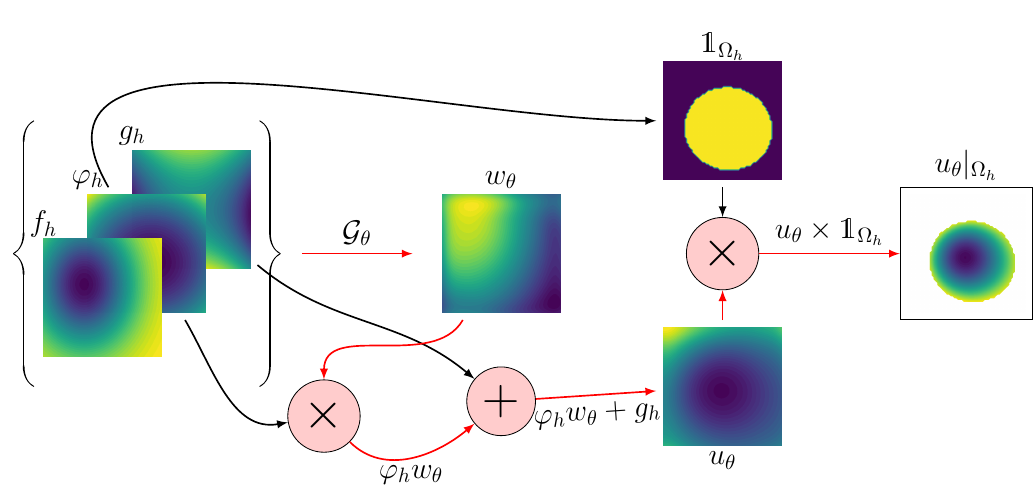}
    \caption{Construction of a prediction of $\varphi$-FEM-FNO to solve \eqref{eq:governing_poisson}.
    }\label{fig:multiplot}
\end{figure}

The goal of the FNO is to construct a parametric mapping
\begin{align*}
    \mathcal{G}_\theta : \quad \bb{R}^{n_x \times n_y \times 3} \quad & \to \quad \bb{R}^{n_x \times n_y \times 1}\,, \\
    (f_h, \varphi_h, g_h) \quad                                       & \mapsto \quad w_\theta\,,
\end{align*}

that approximates the "ground truth" mapping $\mathcal{G}^\dagger$ \eqref{GdaggerVh}.
We predict the $\varphi$-FEM representation $w_h$ with $w_{\theta}$ and
$u_{\theta}=\varphi_h w_{\theta} +g_h$ will be an approximation of $u_h=\varphi_hw_h+g_h$  as described in Figure \ref{fig:multiplot}. Here $\theta$ stands for the numerous parameters that we have to find by minimizing the loss function.

\begin{remark}
    {
        The decision to predict $w_h$ instead of $u_h$ directly stems from the fact that multiplying by  $\varphi_h$ ensures the exact imposition of boundary conditions. Predicting $u_h$ directly would introduce errors at the domain boundaries and could necessitate the inclusion of a loss term specifically targeting the solution's boundary values. In Section \ref{sec:num}, we will demonstrate the impact of predicting $w_h$ versus
        $u_h$ in the first test case, and compare it to a scenario where $u_h$ is predicted, with no loss modification.}
\end{remark}
\subsubsection{The structure of the FNO}

\paragraph{Layer sequence} The mapping $\mathcal G_\theta$ is composed of several sub-mappings, called layers
\begin{equation*}
    \mathcal G_\theta =  N^{-1}\circ Q_{\theta} \circ \mathcal H^{4}_\theta \circ \mathcal H^{3}_\theta \circ \mathcal H^{2}_\theta \circ \mathcal H^{1}_\theta \circ P_{\theta}\circ N \,.
\end{equation*}
Each layer acts on three-dimensional tensors with the third dimension (the number of channels) varying from one layer to another according to the following scheme:
\begin{multline*}
    \mathcal G_\theta:
    \mathbb R^{n_x\times n_y \times 3} \xrightarrow{N }\mathbb R^{n_x\times n_y \times 3} \xrightarrow{P_\theta }\mathbb R^{n_x\times n_y \times n_d} \xrightarrow{\mathcal{H}_\theta^1} \mathbb R^{n_x\times n_y \times n_d} \xrightarrow{\mathcal{H}_\theta^2} \\
    \dots \xrightarrow{\mathcal{H}_\theta^4} \mathbb R^{n_x\times n_y \times n_d}   \xrightarrow{Q_\theta} \mathbb R^{n_x\times n_y \times 1}\xrightarrow{N^{-1}} \mathbb R^{n_x\times n_y \times 1} \,,
\end{multline*}
where $n_d$ is a sufficiently large dimension. A graphic representation of $\mathcal G_\theta$ is given in Fig.~\ref{fig:pipeline}.
The transformations $P_{\theta}$ and $Q_{\theta}$ are respectively an embedding to the high dimensional channel space and a projection to the target dimension, both computed using a neural network (see \cite{paper_FNO}).
The principal components of  FNO are the  4 Fourier layers $\mathcal H^{\ell}_\theta$ having the same structure presented below.

\paragraph{Normalisation $N$ and $N^{-1}$}

To improve the performance of a neural network, it is known that normalizing the inputs and outputs is almost mandatory (see \cite{nastorg:hal-03970501} for example).
So the train data will be standardized channel by channel in the input and unstandardized in the output thanks to the operator $N$ and $N^{-1}$, respectively. This step is explained in detail in Appendix \ref{sec:norm}.

\paragraph{Structure of  the embedding $P_{\theta}$ and the projection $Q_{\theta}$}
The transformation $P_{\theta}$ is made of one fully connected layer of size $n_d$ neurons acting on each node, i.e. for all $i\in\{1,...,n_x\}$, $j\in\{1,...,n_y\}$ and $k\in\{1,...,n_d\}$,

\[
    P_{\theta}(X)_{ijk}= \sum_{k' = 1}^{3}W^{P_\theta}_{kk'}  X_{ijk'} +B^{P_{\theta}}_k\,,
\]

with $W^{P_\theta}\in\mathcal{M}_{n_d,3}(\mathbb{R})$,   $B^{P_{\theta}}\in\mathbb{R}^{n_d}$ some parameters.
The transformation $Q_{\theta}$ is made of two fully connected layers of size $n_Q$ then $1$ acting also on each node, i.e.
$Q_{\theta}=(Q_{\theta,ijk})_{ijk}$ defined for all $X=(X_{ijk})_{ijk}$ by, for all $i\in\{1,...,n_x\}$, $j\in\{1,...,n_y\}$,

\[
    Q_{\theta}(X)_{ij}=\left[\sum_{k=1}^{n_Q} W^{Q_{\theta,2}}_{1k}\sigma\left(\sum_{k' = 1}^{n_d}W^{Q_{\theta,1}}_{kk'}  X_{ijk'} +B^{Q_{\theta,1}}_k\right)\right] + B^{Q_{\theta,2}}\,,
\]
with $W^{Q_{\theta,1}}\in\mathcal{M}_{n_d,n_Q}(\mathbb{R})$,   $B^{Q_{\theta,1}}\in\mathbb{R}^{n_Q}$, $W^{Q_{\theta,2}}\in\mathcal{M}_{n_Q,1}(\mathbb{R})$,   $B^{Q_{\theta,2}}\in\mathbb{R}$ some parameters and $\sigma$ is an activation function applied term by term. We choose the GELU (Gaussian Error Linear Unit) function given by  $f(x) = x \varphi(x)$ with $\varphi(x) = P(X\leqslant x)$ where $X \sim \mathcal{N} (0,1)$,
as in the original implementation of the FNO\footnote{\url{https://github.com/neuraloperator/neuraloperator}} and of the Geo-FNO\footnote{\url{https://github.com/neuraloperator/Geo-FNO}}.

\paragraph{Structure of  Fourier layers $\mathcal H^{\ell}_\theta$}
A layer $\mathcal H^{\ell}_\theta$ is made of two sub-layers organized as follows (see \cite{paper_FNO}):
\[
    \mathcal H^{\ell}_\theta (X) =  \sigma \big( \mathcal C^{\ell}_\theta (X) + \mathcal B^{\ell}_\theta (X) \big)\,,
\]
where
\begin{itemize}

    \item $\mathcal C^{\ell}_\theta$ is a layer defined  by
          \begin{equation*}
              \mathcal C^{\ell}_{\theta}(X) = \mathcal F^{-1}\Big(W^{\mathcal{C}^{\ell}_{\theta}}  \mathcal F(X)\Big) \in  \mathbb{R}^{n_x\times n_y\times n_d \times n_d}\,,
          \end{equation*}
          with $W^{\mathcal{C}^{\ell}_{\theta}} \in \bb{C}^{n_x\times n_y\times n_d\times n_d
              }$
          some parameters
          and $\mathcal F$, $\mathcal F^{-1}$ stand for the real Fast Fourier Transform (RFFT) and its inverse:
          for all $i\in\{1,...,n_x\}$, $j\in\{1,...,n_y\}$ and $k\in\{1,...,n_d\}$,
          \[
              \mathcal F(X)_{ijk} = \sum_{i'j'} X_{i'j'k} e^{2\sqrt{-1}\pi \left( {\frac{ii'}{n_x}} + {\frac{jj'}{n_y}}\right)}
              \,,
          \]
          and for $Y\in\mathbb{C}^{n_x\times n_y\times n_d}$
          \[
              \mathcal F^{-1}(Y)_{ijk} =
              \sum_{i'j'} Y_{i'j'k} e^{-2\sqrt{-1}\pi  \left( {\frac{ii'}{n_x}} + {\frac{jj'}{n_y}}\right)}
              \,.
          \]
    \item $\mathcal B^{\ell}_\theta=(\mathcal B^{\ell}_{\theta,ijk})_{ijk}$ is the "bias-layer" defined for all $X=(X_{ijk})_{ijk}$ by, for all $i\in\{1,...,n_x\}$, $j\in\{1,...,n_y\}$ and $k\in\{1,...,n_d\}$,

          \[
              \mathcal{B}^\ell_{\theta}(X)_{ijk}= \sum_{k' = 1}^{n_d}W^{\mathcal{B}^\ell_\theta}_{kk'}  X_{ijk'} +B^{\mathcal{B}^\ell_{\theta}}_k\,,
          \]

          with $W^{\mathcal{B}^\ell_\theta} \in\mathcal{M}_{n_d}(\mathbb{R})$ and $B^{\mathcal{B}^\ell_{\theta}}\in\mathbb{R}^{n_d}$.
\end{itemize}

The coefficients {of} $W^{\boldsymbol{\cdot}}$ and $B^{\boldsymbol{\cdot}}$ are the trainable parameters, which have some constraints given below.

\paragraph{Constraint on the parameters}

\begin{itemize}
    \item \textbf{Symmetry}:
          To obtain a real matrix $\mathcal C^{\ell}_{\theta}(X)$, we impose to  $W^{\mathcal{C}^{\ell}_{\theta}}$ the Hermitian symmetry, i.e. $W^{\mathcal{C}^{\ell}_{\theta}}_{n_x-i,n_y-j,k}=\overline{W}^{\mathcal{C}^{\ell}_{\theta}}_{i,j,k}$.
          In practice, since we use a specific implementation of the FFT called RFFT (Real-FFT): the discrete Fourier coefficients are stored in arrays of size $n_x\times (n_y/2+1)$ (integer division) and are automatically symmetrized when performing the inverse transformation. Hence, there is no precaution to take for the matrix $\hat W$ in practice. Simply it must be of size $n_x\times (n_y/2+1)$.

    \item \textbf{Low pass filter}: The solutions to our problem \eqref{eq:governing_poisson} are usually rather smooth. Hence, when we perform the RFFT, the "very high" frequencies can be neglected. They simply participate in the fact that the RFFT is a bijection. Typically, a good approximated solution can be recovered, keeping only the $m \times m$ first "low" frequencies.
\end{itemize}

\begin{remark}
    [Number of parameters]
    An interesting aspect of the FNO is the reasonable number of parameters to optimize. Indeed, since we truncate the high frequencies, for each $\mathcal{C}^l_\theta$ the number of parameters to optimize is less than $n_x \times n_y \times n_d \times n_d$. In fact, the total number of parameters $n_\theta$ does not depend on the resolution and is given by
    \[
        n_\theta = \overbrace{4 \times n_d}^{P_\theta \ : \ 3\times n_d + n_d} + \  4 \times \underbrace{(\overbrace{2 \times n_d^2 \times m^2}^{\substack{\text{$\mathcal{C}^l_\theta$ : by truncation of} \\ \text{high frequencies}}} + \overbrace{n_d^2 + n_d}^{\mathcal{B}^l_\theta} )}_{\mathcal{H}^l_\theta} + \overbrace{(n_d + 2) \times n_Q + 1}^{Q_\theta \ : \ n_d \times n_Q + n_Q + n_Q  \times 1 + 1}\,. \]
    In the following first test case, for the chosen hyperparameters this will represent 324577 parameters to optimize.
\end{remark}

\begin{remark}
    Once trained, FNOs can work with input data of arbitrary dimensions $n_x$, $n_y$. This property of multi resolution is due to a special structure of FNO.  In fact, FNO is presented in \cite{paper_FNO,Neural_operator} as an approximation of mappings acting on infinite dimensional function spaces. It is thus not surprising that its discretization can be done, in principle, on any mesh.
\end{remark}

\paragraph{Border issues}
An issue of the RFFT applied to non-periodic functions is the Gibbs phenomenon: some oscillations appear on the border.
To erase them, we apply a padding technique: we extend the matrices, adding entries all around, before performing the convolution. After the convolution, we restrict the matrix to its original shape to partially erase the oscillations (see the PyTorch documentation\footnote{\url{https://pytorch.org/docs/stable/generated/torch.nn.ReflectionPad2d.html}} for an example).

\section{The proposed architecture}\label{sec:main idea}

\begin{figure}
    \centering
    \includegraphics[width=0.8\textwidth]{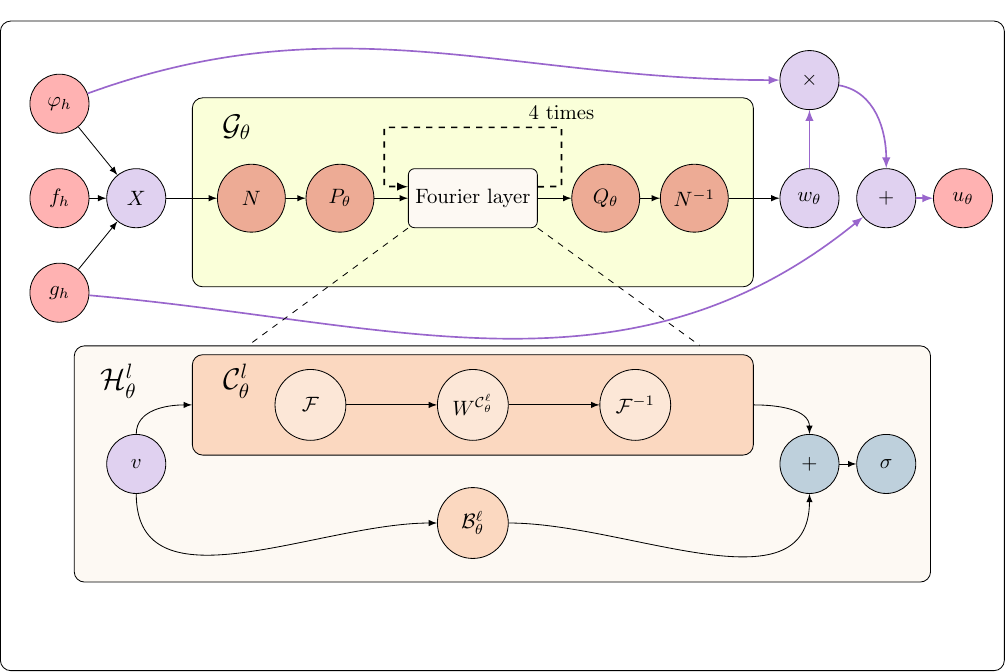}
    \caption{The $\varphi$-FEM-FNO pipeline to solve \eqref{eq:governing_poisson}. Illustration based on the representation of \cite{paper_FNO}. The upper part represents the entire pipeline, and the lower part is a zoom on a Fourier layer. The red circles correspond to the inputs provided by the user and the output returned by our $\varphi$-FEM-FNO. We represent the inputs and outputs seen by the FNO in purple, where $X = (f_h,\varphi_h, g_h)$. In orange, $P_\theta$ and $Q_\theta$ are two transformations parameterized by neural networks. Moreover, $\mathcal{F}$ and $\mathcal{F}^{-1}$ are respectively the Fourier and inverse Fourier transforms. In blue, $\sigma$ is the activation function. Finally, black arrows correspond to steps inside our FNO, and purple arrows to steps outside the FNO.}\label{fig:pipeline}
\end{figure}

As previously said, the main objective of this paper is to present a novel approach that combines two existing methods: $\varphi$-FEM and the Fourier Neural Operator (FNO). This method aims to leverage the high precision of $\varphi$-FEM while utilizing the FNO's ability to generate nearly instantaneous predictions after training. Such a combination enables the approach to be employed effectively in real-time simulations. An overview of the entire pipeline is provided in Fig.~\ref{fig:pipeline}, where the input data and the final output are highlighted in red.

\begin{figure}
    \centering
    \includegraphics[width=\textwidth]{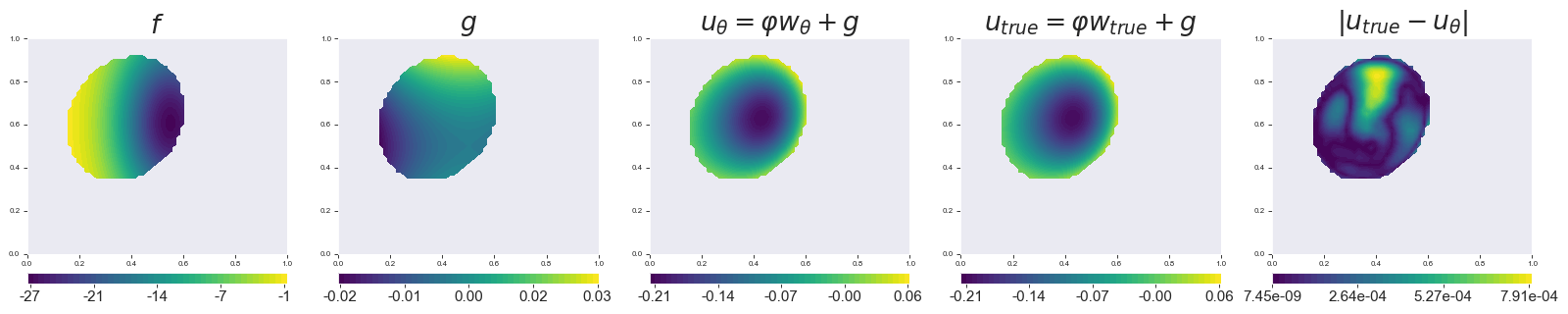}
    \caption{\textbf{Test case 1.} Example of result among the validation sample with an error in the norm \eqref{eq:norm_L2_tensors} of $2.5\times 10^{-3}$, corresponding to the median error. }
    \label{fig:output_fno}
\end{figure}

By construction, a prediction of the FNO will be given on the same cartesian regular grid of the inputs. Since we are interested in the solution only over $\Omega_h$, we need to define a loss function acting only on the corresponding pixels. An example of data and truncated output of our approach is represented in Fig.~\ref{fig:output_fno}.

\paragraph*{Loss function.}

Let $N_{\text{data}}$ be the size of a considered sample of data. We denote $U_{\text{true}} = (u_{\text{true}}^n)_{n = 0, \dots, N_{\text{data}}}$ where $u_{\text{true}}^n = \varphi_h^n w_h^n + g_h^n$, the ground truth solution and $U_{\theta} = (u_{\theta}^n)_{n = 0, \dots, N_{\text{data}}}$ with $u^n_\theta = \varphi_h^n \mathcal{G}_\theta(f_h^n, \varphi_h^n, g_h^n) + g_h^n = \varphi_h^n w_\theta^n + g_h^n$ the output of $\varphi$-FEM-FNO.

\begin{figure}
    \centering

    \begin{tikzpicture}[scale=0.7, every node/.style={minimum size=0.8cm, draw, anchor=center}]

        \foreach \dx/\dy in {2/1,3/1,4/1,5/1,
                1/2,2/2,3/2,4/2,5/2,6/2,
                1/3,2/3,3/3,4/3,5/3,6/3,
                1/4,2/4,3/4,4/4,5/4,6/4,
                1/5,2/5,3/5,4/5,5/5,6/5,
                2/6,3/6,4/6,5/6
            } {
                \fill[gray!50] (\dx,\dy) rectangle ++(1,1);
            }

        \foreach \dx/\dy in {3/2,4/2,
                2/3,3/3,4/3,5/3,
                2/4,3/4,4/4,5/4,
                3/5,4/5} {
                \fill[blue!80] (\dx,\dy) rectangle ++(1,1);
            }
        \foreach \x in {0,1,2,3,4,5,6,7,8} {
                \draw[black] (\x,0) -- (\x,8);
                \draw[black] (0,\x) -- (8,\x);
            }
        \draw[red] (4,4) circle[radius=75pt];
    \end{tikzpicture}
    \caption{In red, the real boundary of an example domain.
        In blue and gray, the set $\mathcal{S}_0$. In gray, $\mathcal{S}_1$. }\label{fig:masks}
\end{figure}
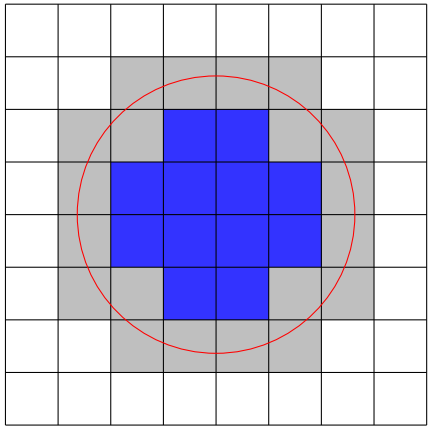

The loss to be optimized is an approximation of the average $H^1$ error, given by

\begin{equation}\label{eq:loss}
    \mathcal{L}\left(U_{\text{true}}; U_{\theta}\right) = \frac{1}{N_\text{data}}\sum_{n=0}^{N_\text{data}}
    \left( \mathcal{E}_0( u^n_{\text{true}} ; u^n_{\theta}) + \mathcal{E}_1( u^n_{\text{true}} ; u^n_{\theta})\right)\,,
\end{equation}
where
\begin{align*}
    \mathcal{E}_0(u^n_{\text{true}}; u^n_{\theta}) & =
    \| u^n_{\text{true}} - u^n_{\theta} \|^2_{0, \mathcal{S}_0^n}
    \,,                                                \\
    \intertext{and}
    \mathcal{E}_1(u^n_{\text{true}}; u^n_{\theta}) & =
    \| \nabla_x^h u^n_{\text{true}} - \nabla_x^h u^n_{\theta} \|^2_{0,\mathcal{S}_1^n}
    + \| \nabla_y^h u^n_{\text{true}} - \nabla_y^h u^n_{\theta} \|^2_{0, \mathcal{S}_1^n} \,,
\end{align*}
where $\nabla^h$ is the centered finite difference approximation of the gradient and $\mathcal{S}_0$ is the set of pixels corresponding to the vertices of $\Omega_h$, and $\mathcal{S}_1$ is the set of pixels of $\mathcal{S}_0$ deprived of a layer of pixels (constructed using the 8-th neighborhood, see Fig.~\ref{fig:masks}).

\begin{remark}
    In \eqref{eq:loss}, we compute the loss with respect to $u_\text{true}^n$ and not $w_\text{true}^n$. This way of computing the loss does not mean that our FNO will predict $u_\theta^n$. It only means that we will predict a solution $w_\theta$ such that, multiplied by input function $\varphi_h$ and added to $g_h$, the result will be close to $u_\text{true}^n$. In the following test cases, we will illustrate numerically the difference between using $w_\text{true}$ and $u_\text{true}$ in the loss.
\end{remark}

\section{Numerical results} \label{sec:num}

Let us now illustrate the efficiency of our technique by numerical test cases. We will first solve the Poisson equation \eqref{eq:governing_poisson} considering the case of parametric domains, using varying elliptic domains to illustrate the accuracy and the fastness of the method compared to five other methods. Then, we will extend our study to more complex shapes and finally to the resolution of a hyperelastic problem \eqref{eq:governing_hyperelast}.

We will fix $n_d=20$ (number of neurons acting on each node), $n_Q=128$ (number of neurons in the first layer of the projection $Q_{\theta}$), and $m=10$ (the number of low frequencies considered in the low pass filter).
In the data produced by $\varphi$-FEM, the parameter $\sigma_D$ is fixed to $1$ for the two first test cases.

\paragraph{Implementation details}
All the simulations were executed on a laptop with an \texttt{Intel Core i7-12700H CPU}, $32$Gb of memory, and an \texttt{NVIDIA RTX A2000 GPU} with $8$Gb of memory. The data were generated using the python finite element library {DOLFINx (\cite{DOLFINx,basix1,basix2,uflX})} and the FNO is implemented\footnote{All the codes and datasets used in this paper are available at \url{https://github.com/KVuillemot/PhiFEM_and_FNO}.} using the \textit{Pytorch}\cite{pytorch} library. Moreover, we will use an ADAM optimizer with an initial learning rate $\alpha=0.0005$, $\beta_1 = 0.9$, $\beta_2 = 0.999$, and $\varepsilon = 10^{-7}$ to train the operator (see Appendix \ref{Ap:algo}  Algo. \ref{algo:ADAM}). During training, the learning rate is reduced when the loss on the validation sample does not decrease over several epochs.
The algorithm of the training loop is presented in Appendix \ref{Ap:algo} Algo. \ref{algo:training}.

\paragraph*{Evaluation metrics.}
To evaluate the performance of the FNO, we define two different metrics, allowing us to compute the $L^2$ relative errors:
\begin{itemize}
    \item The first metric that will be used to compute the error between a FNO solution and a ground truth solution, is given by
          \begin{equation}\label{eq:norm_L2_tensors}
              {E_1(u_{\text{true}}, u_\theta) := }  \sqrt{\frac{\mathcal{E}_0(u_{\text{true}};u_{\theta})}{\mathcal{N}_0(u_{\text{true}})}}\,,
          \end{equation}
          where $u_{\theta} = \varphi_h \mathcal{G}_\theta (\varphi_h, f_h, g_h) + g_h $ and $u_\text{true} = \varphi_h w_h + g_h$. Moreover, we denote $\mathcal{L}_0(\cdot)$ the average value of this metric among a given dataset (train, validation, etc).

    \item The second metric will be used to compute the errors with respect to fine standard finite element solutions $u_\text{ref}$ and is defined by

          \begin{equation}\label{eq:norm_L2_fenics}
              E_2(u_{\text{ref}}, u_\theta) := \frac{\|\Pi_{\Omega_{\text{ref}}} u_{\theta} - u_{\text{ref}} \|_{0, \Omega_{\text{ref}}}}{\|u_{\text{ref}} \|_{0, \Omega_{\text{ref}}}} =
              \sqrt{\frac{\int_{\Omega_{\text{ref}}} (\Pi_{\Omega_{\text{ref}}} u_{\theta}-u_{\text{ref}})^2 \ \mathrm{d}x}{\int_{\Omega_{\text{ref}}} u_{\text{ref}}^2 \ \mathrm{d} x}} \,,
          \end{equation}

          where $\Pi_{\Omega_{\text{ref}}}$ denotes an approximation of the $L^2$-orthogonal projection on the reference domain $\Omega_{\text{ref}}$ (fine conformed mesh of $\Omega$).

\end{itemize}

\subsection{The Poisson-Dirichlet equation on varying ellipses}

Let us first consider the simple case of the Poisson equation \eqref{eq:governing_poisson} on elliptic domains given by the level-set functions

\begin{multline}\label{eq:level_set_moving} \varphi_{(x_0, y_0, l_x, l_y, \theta)}(x,y) = -1 +
    \frac{((x-x_0)\cos(\theta) + (y-y_0)\sin(\theta))^2}{l_x^2}  \\
    + \frac{((x-x_0)\sin(\theta) - (y-y_0)\cos(\theta))^2}{l_y^2}  \,, \end{multline}
with
\[ x_0,\  y_0 \sim \mathcal{U}([0.2,0.8])\,, \quad l_x, \ l_y \sim \mathcal{U}([0.2,0.45]) \quad \text{ and } \theta \sim \mathcal{U}([0, \pi])\,.
\]

The equation \eqref{eq:level_set_moving} defines an ellipse centered in $(x_0, y_0)$ of semi-major axis $l_x$ and semi-minor axis $l_y$, rotated by an angle $\theta$ around the center of the ellipse, as illustrated for example in Fig.~\ref{fig:domain_submeshes} (left).
We apply a rejection sampling method on the previous parameters to ensure that each domain is entirely lying within the unit square.
The random functions $f$ and $g$ of \eqref{eq:governing_poisson} are given by
\begin{equation}
    \label{eq:random_f}
    f_{(A,\mu_0,\mu_1, \sigma_x, \sigma_y)}(x,y) = A \exp\left(- \frac{(x-\mu_0)^2}{2\sigma_x^2} - \frac{(y-\mu_1)^2}{2\sigma_y^2}\right)\,,
\end{equation}
and
\begin{equation}
    \label{eq:random_g}
    g_{(\alpha, \beta)}(x,y) = \alpha \left((x-0.5)^2 - (y-0.5)^2\right) \cos\left( \beta y \pi \right)\,,
\end{equation}
where $A\sim \mathcal{U}([-30, -20] \cup [20,30])$, $(\mu_0,\mu_1) \sim \mathcal{U}([0.2,0.8]^2\cap\{ \varphi < - 0.15 \})$, $\sigma_x$, $\sigma_y \sim \mathcal{U}([0.15,0.45])$ and $\alpha$, $\beta \sim \mathcal{U}([-0.8,0.8])$.

We generate a set of data of size 2100, split into a training set of size 1500, a validation set of size 300, and a test set of size 300. The training set is then divided into batches of size $32$ (number of data considered in the loss function in one computation of the gradient) at each of the $2000$ epochs
(number of loops over all the batches), as explained in Algorithm~\ref{algo:training}.

\begin{remark}[Data generation]
    Note that to generate the data we use $\bb{P}^2$ interpolations of $f$ and $\varphi$, considering that we can use a maximal information. However, to compare the methods, since the FNO approaches are based on the nodal values, we will use only the nodal values of the functions for the FEMs-based methods, to have comparable results.
\end{remark}

\begin{figure}
    \centering
    \includegraphics[width=0.9\textwidth]{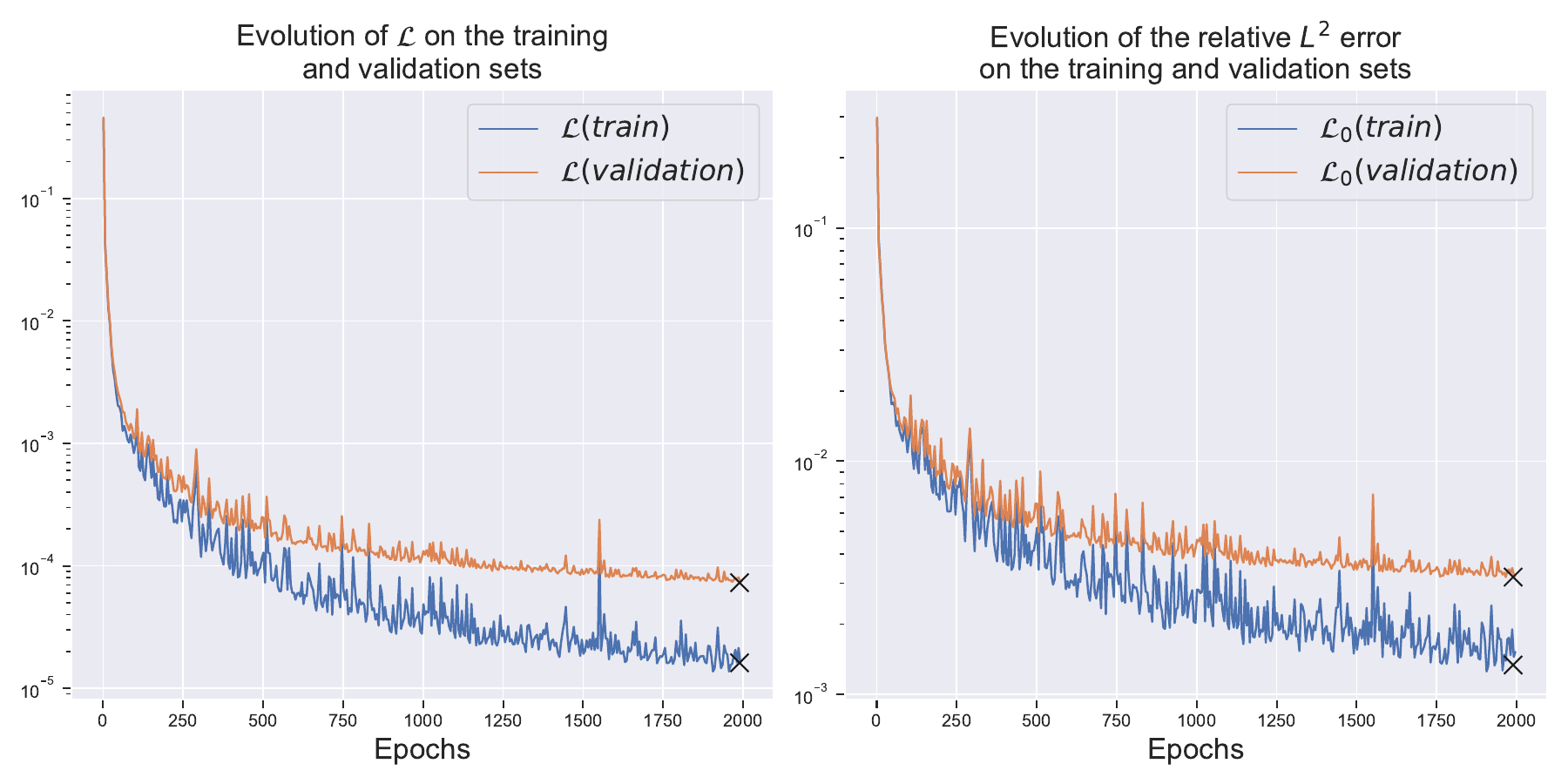}
    \caption{\textbf{Test case 1.} On the left-hand side (resp. right-hand side), we represent the evolution of the cost function $\mathcal{L}$ (resp. the relative $L^2$ error) on a subset of the training set and on the validation set.}
    \label{fig:loss_test_case_1}
\end{figure}

\paragraph*{Results on the validation sample}
We represent in Fig.~\ref{fig:loss_test_case_1} (left) the evolution of the loss function $\mathcal{L}$ on a random subset of the training dataset and on the validation dataset, both of size 300. In Fig.~\ref{fig:loss_test_case_1} (right), we represent the evolution of $\mathcal{L}_0$ on the same samples of data.

During the training, we select the model minimizing the loss function $\mathcal{L}$ on the validation set. This model will be considered to be the optimal model returned by the training and will be used in the third part of this test case.

\begin{figure}
    \centering
    \includegraphics[width=0.8\textwidth]{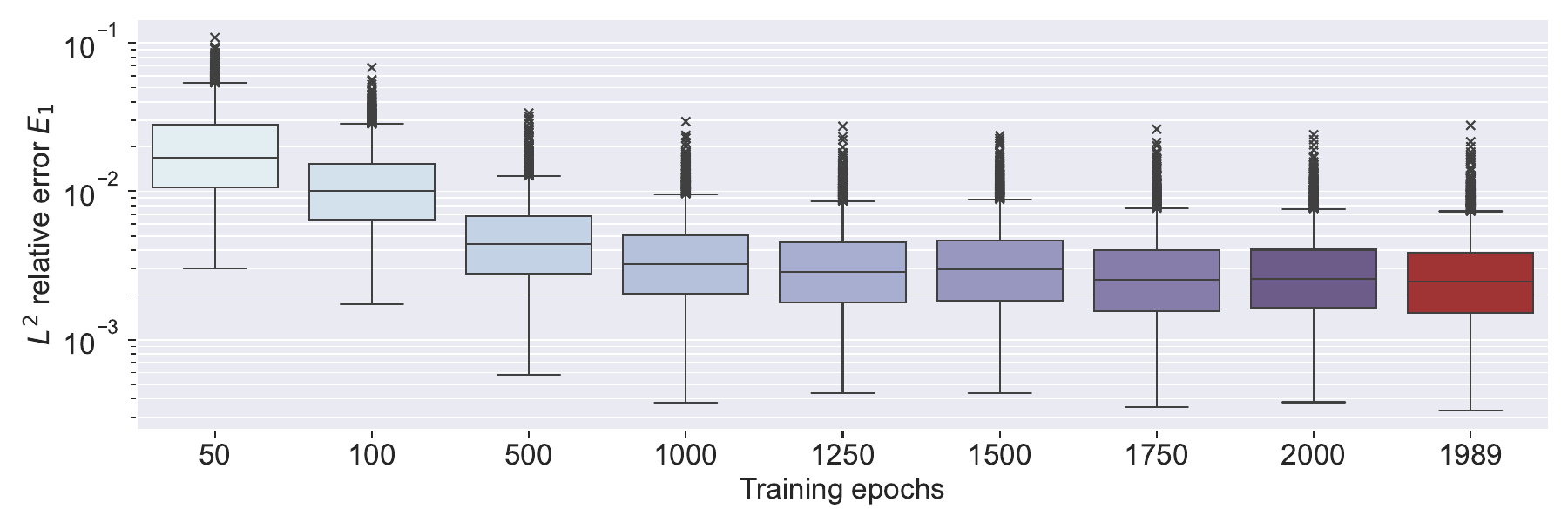}
    \caption{\textbf{Test case 1.}
        Evolution of the \eqref{eq:norm_L2_tensors} errors on 2500 test data at different steps of the training. The optimal model is represented in red.}
    \label{fig:test_case_1_boxplot_new_data_phifem}
\end{figure}

\paragraph*{Validation of the model of a first test dataset.}
We now address a second crucial aspect: evaluating the error of the models on a test dataset, using the norm \eqref{eq:norm_L2_tensors}. This evaluation ensures that the operator is well trained and performs consistently on new data, behaving as on the validation data. Additionally, it serves as a final validation of the optimality of the selected best model. To achieve this, we consider 2500 new data and compute the error in the norm \eqref{eq:norm_L2_tensors} at several steps of the training, including the optimal one. The results, shown in Fig.~\ref{fig:test_case_1_boxplot_new_data_phifem}, illustrate that the selected optimal model is the best among those tested.

\paragraph*{{Comparison of $\varphi$-FEM-FNO, FEM's based methods and other ML-based methods }}

We now turn to the key numerical results to illustrate the advantages of our $\varphi$-FEM-FNO technique.
In this part, we compare $\varphi$-FEM-FNO with several other approaches, to highlight the effectiveness of our method:
\begin{itemize}
    \item $\varphi$-FEM-FNO: we call the previous optimal model with 1500 training data of resolution $64 \times 64$ (corresponding to a cell size $h\approx 0.022$);
    \item $\varphi$-FEM-FNO 2: we apply the same process as for $\varphi$-FEM-FNO, but predicting directly $u_\theta$ instead of $w_\theta$, i.e. we define a new operator
          \begin{align*}
              \mathcal{G}_\theta : \quad \bb{R}^{n_x \times n_y \times 3} \quad & \to \quad \bb{R}^{n_x \times n_y \times 1}\,, \\
              (f_h, \varphi_h, g_h) \quad                                       & \mapsto \quad u_\theta\,.
          \end{align*}
          We use the same loss function (i.e. $\mathcal{L}$ defined in \eqref{eq:loss}) except that $u^n_\theta$ is directly the prediction;
    \item $\varphi$-FEM-UNET: we have adapted the previous framework to another well-known Neural Network architecture, namely the UNet architecture (see \cite{ronneberger2015u}). To train this network, we use the loss function $\mathcal{L}$ defined in \eqref{eq:loss}. Note that this model represents much more parameters to optimize than the FNO-based methods. For this test case, we optimize a total of 7753025 parameters (20 times more than for $\varphi$-FEM-FNO).
    \item $\varphi$-FEM: we apply the operator $\mathcal{G}^\dagger$, with background meshes of resolution $64 \times 64$ (corresponding to a cell size $h\approx 0.022$), taking $\sigma_D =1$ and $\bb{P}^1$ finite elements;
    \item Standard FEM: we use a standard $\bb{P}^1$ finite element method to solve the problems on meshes with cell size $h \approx 0.022$, corresponding to the resolution used for the other approaches;
    \item Standard-FEM-FNO: we use a FNO trained with standard $\bb{P}^1$ FEM solutions on meshes of sizes $h\approx 0.022$, interpolated on Cartesian grids of size $64 \times 64$ as data. The loss function used to train Standard-FEM-FNO is the relative $H^1$ norm and the operator is trained during 2000 epochs.
    \item Geo-FNO: we have trained a Geo-FNO, adapting the approach of \cite{FNO_general_geometries} (see the implementation on the GitHub\footnote{\url{https://github.com/neuraloperator/Geo-FNO/blob/main/elasticity/elas_geofno_v2.py}}) to match our test case, using as input of the operator a set of 2600 points and the values of $f$ and $g$ at each of these points. We used 2600 points to obtain an average cell size close to $0.02$. The operator has been trained during 2000 epochs, using the $L^2$ relative norm.

\end{itemize}

The five ML-based methods have been trained using the same dataset, adapted to each method (i.e., during the data generation steps, we consider the same set of parameters, and the same hyperparameters during training.)

\begin{remark}[Implementation aspect.] To compare our method with a standard finite element method, we need to construct conforming meshes using the values of the level-set function $\varphi$. The creation of a mesh using a level-set function is not directly possible with DolfinX. Thus, we need to create such meshes manually. For this step, we use \textit{pymedit}\footnote{\url{https://pypi.org/project/pymedit/}} with the Mmg platform\footnote{\url{https://www.mmgtools.org/}}. We refer the reader to the GitHub repository\footnote{\url{https://github.com/KVuillemot/PhiFEM_and_FNO/blob/main/install_and_use_mmg.md}} for details on the installation and two examples of use.
\end{remark}

\begin{figure}
    \centering
    \includegraphics[width=0.9\textwidth]{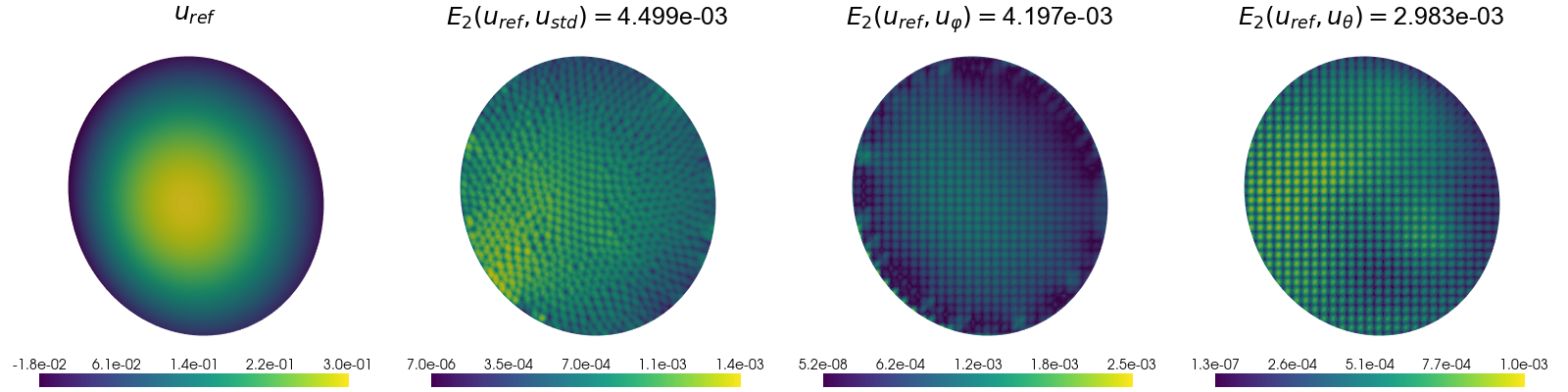}
    \caption{\textbf{Test case 1.} Reference solution ($u_{\text{ref}}$).
        Outputs of standard FEM ($u_{\text{std}}$) , $\varphi$-FEM ($u_{\varphi}$), and ($\varphi$-FEM-FNO $u_{\theta}$). The presented test case corresponds to the median one among the FNO $L^2$ relative errors.}
    \label{fig:test_case_1_outputs}
\end{figure}

\begin{figure}
    \centering
    \includegraphics[width=0.495\textwidth]{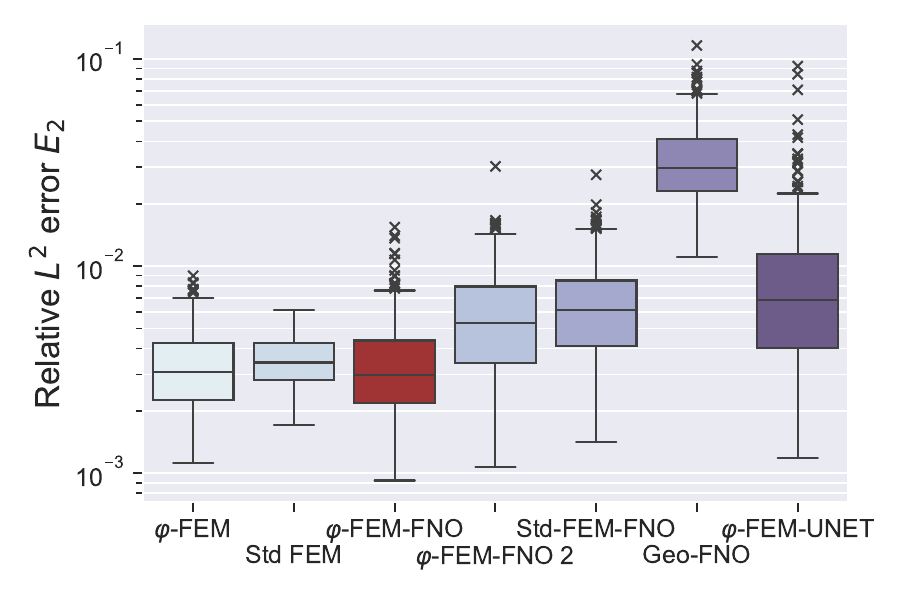}
    \includegraphics[width=0.495\textwidth]{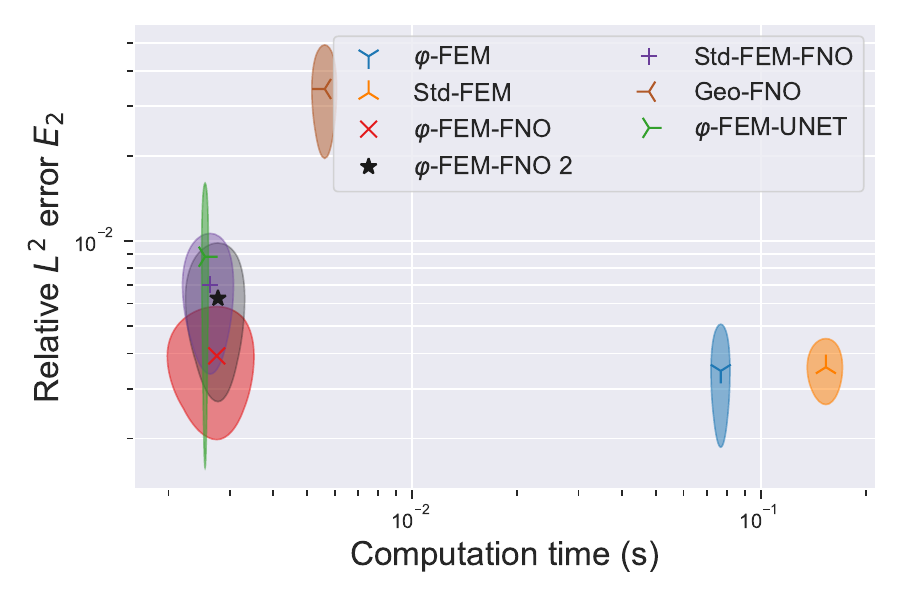}
    \caption{\textbf{Test case 1.} Left: Relative $L^2$ errors of the methods. Right: Relative $L^2$ errors, with respect to the computation times.}\label{fig:compare_methods_test_case_1}
\end{figure}

To compare $\varphi$-FEM-FNO with other methods, we evaluate the best models of each machine learning-based approach using a test sample of size 300. The predicted solutions are projected onto a reference fine mesh with cell sizes of approximately $h_\text{ref} \approx0.005$, as shown in Fig.~\ref{fig:test_case_1_outputs}. Errors are computed using the norm defined in \eqref{eq:norm_L2_fenics}, with a fine standard finite element solution serving as the reference. Fig.~\ref{fig:compare_methods_test_case_1} (left) demonstrates that the trained
$\varphi$-FEM-FNO achieves a precision comparable to FEM-based methods. Furthermore, $\varphi$-FEM-FNO is approximately twice as precise as Standard-FEM-FNO and ten times more precise than Geo-FNO. Additionally, $\varphi$-FEM-FNO outperforms $\varphi$-FEM-UNET, highlighting the advantages of the {FNO}
over UNet architectures.

Finally, while $\varphi$-FEM-FNO-2 also performs better than Standard-FEM-FNO and Geo-FNO, its precision is slightly lower than that of $\varphi$-FEM-FNO.

In Fig.~\ref{fig:compare_methods_test_case_1} (right), each marker represents the average error of a method plotted against the average computation time (in seconds). The shaded regions are constructed using the standard deviation of computation times (x-axis) and relative errors (y-axis) to indicate variability.

For the $\varphi$-FEM method, the computation time includes the sum of the following components: selecting and constructing $\Omega_h$ and $\Omega_h^\Gamma$ (including background mesh construction), interpolation times for $f$, $\varphi$, and $g$, assembling the finite element matrix, and solving the linear system. For standard FEM, the computation time accounts for mesh construction, interpolation times for $f$ and $g$, matrix assembly, and solution of the linear system. In contrast, for ML-based methods, the computation time is the model inference time.

The results clearly demonstrate that ML-based methods are significantly faster than FEM-based approaches. Specifically, they highlight that while $\varphi$-FEM-FNO achieves nearly the same precision as FEM-based methods, it is approximately 100 times faster.

\subsection{Second test case: Poisson equation on random complex shapes}

\begin{figure}
    \centering
    \includegraphics[width=0.49\textwidth]{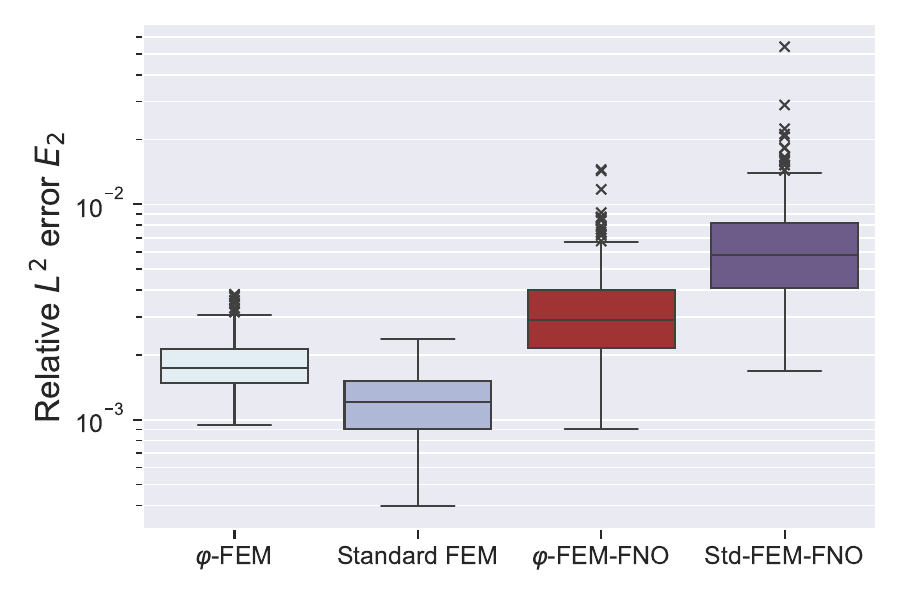}
    \includegraphics[width=0.49\textwidth]{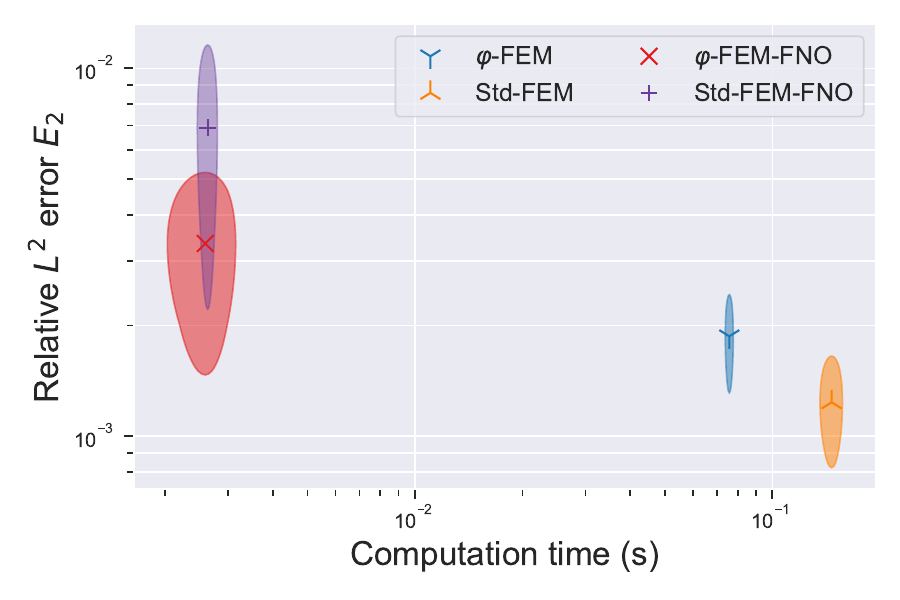}
    \caption{\textbf{Test case 2.} Left: comparison of the three methods on 300 new data. Right: relative $L^2$ error against computation time.}
    \label{fig:results_random_shapes}
\end{figure}

\begin{figure}
    \includegraphics[width=0.49\textwidth]{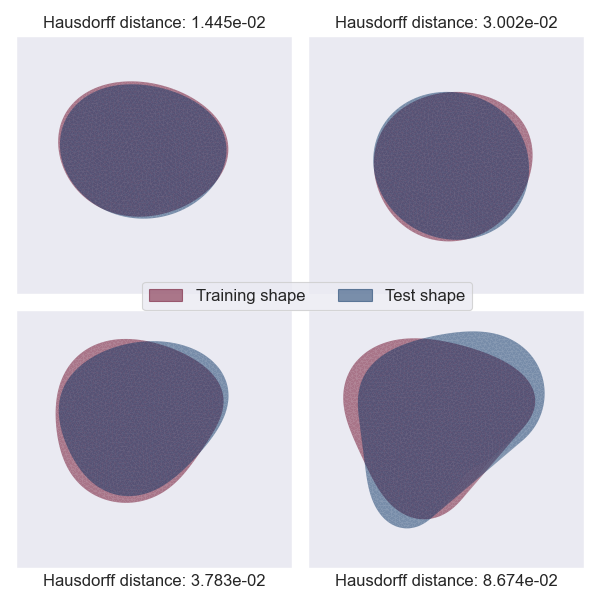}
    \includegraphics[width=0.49\textwidth]{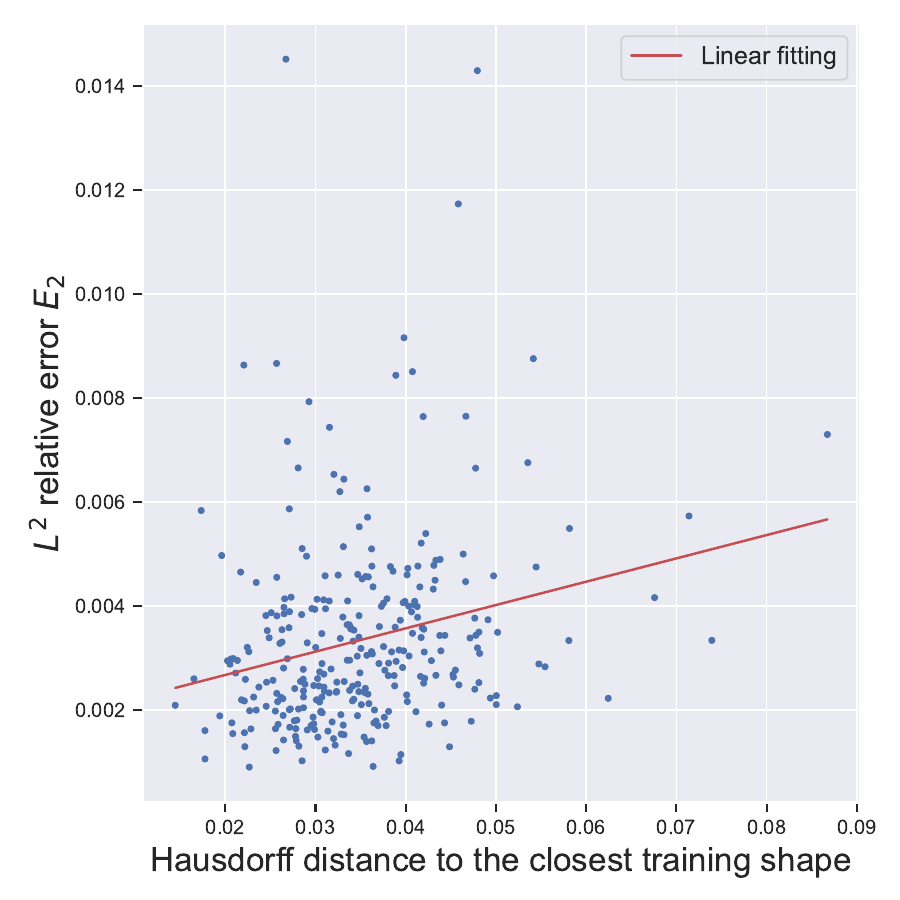}
    \caption{\textbf{Test case 2.} Left: Examples of test geometries. Each represented training shape is the closest one to the considered test shape, in the sense of Hausdorff distance. Right: $L^2$ relative errors against the Hausdorff distance to the closest training shape.}
    \label{fig:results_hausdorff_random_shapes}
\end{figure}

{We now consider a more complex test case while still solving \eqref{eq:governing_poisson}. The functions $f$ and $g$ are defined as in \eqref{eq:random_f} and \eqref{eq:random_g}, with $f$ restricted to positive values. This time, however, we explore more intricate and diverse random shapes. To generate these geometries, we use random level-set functions constructed as a sum of three Gaussian functions. The level-set functions $\varphi$ are defined as:
    \begin{equation}\label{eq:random_phi_gaussian}
        \varphi(x,y) = - \psi(x,y) + 0.5 \max_{(x,y) \in [0,1]^2} \psi(x,y)\,,
    \end{equation}
    with
    \[
        \psi(x,y) = \sum_{k=1}^3 \exp\left( - \frac{(x-x_k)^2}{2 \sigma_k} - \frac{(y-y_k)^2}{2 \gamma_k} \right)\,,
    \]
    where the parameters $x_k$, $y_k$, $\sigma_k$ and $\gamma_k$ are sampled using a Latin Hypercube \cite{latin_hypercube}, along with the parameters of the functions $f$ and $g$.

    The training hyperparameters remain the same as in the first test case, except for the batch size, which is set to 8.

    \begin{figure}
        \centering
        \includegraphics[width=\textwidth]{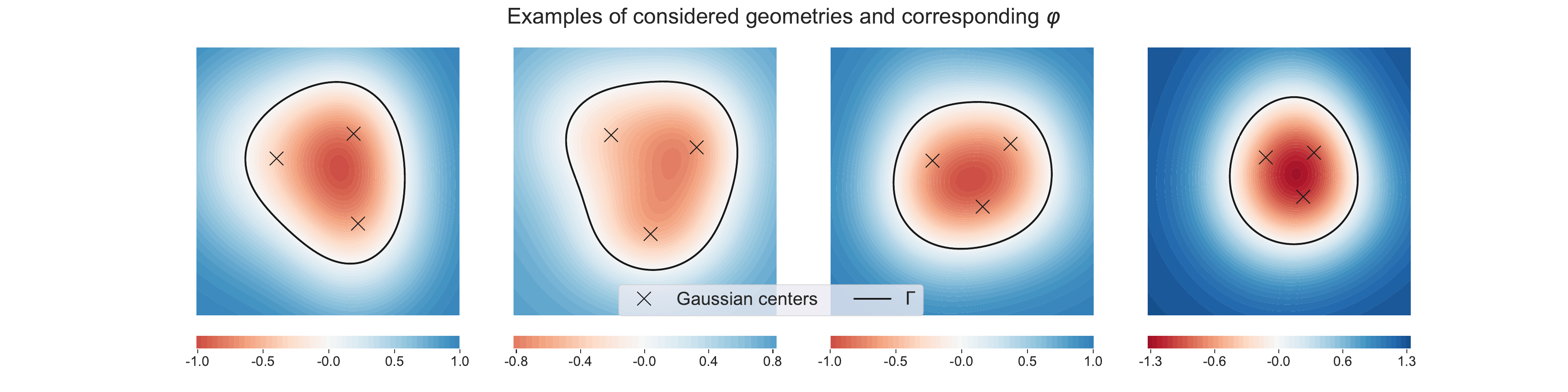}
        \caption{\textbf{Test case 2.} Examples of considered level-set functions and domains, using \eqref{eq:random_phi_gaussian}. The black crosses mark the centers of the Gaussian functions.}
        \label{fig:random_domains}
    \end{figure}

    We represent 4 examples of such level-set functions and given domains, in Fig.~\ref{fig:random_domains}.
    We train the operator for 2000 epochs using 500 training data and 300 validation data. As in the first test case, we use the average $H^1$ norm \eqref{eq:loss} as loss function to minimize. We evaluate the performance of the model on 300 test data, comparing the method to a standard finite element approach, to $\varphi$-FEM, and Standard-FEM-FNO. To evaluate the performance of the methods, as in the first test case, we use a reference standard fine solution, as depicted in Fig.~\ref{fig:output_random_shapes}.

    \begin{figure}
        \centering
        \includegraphics[width=\textwidth]{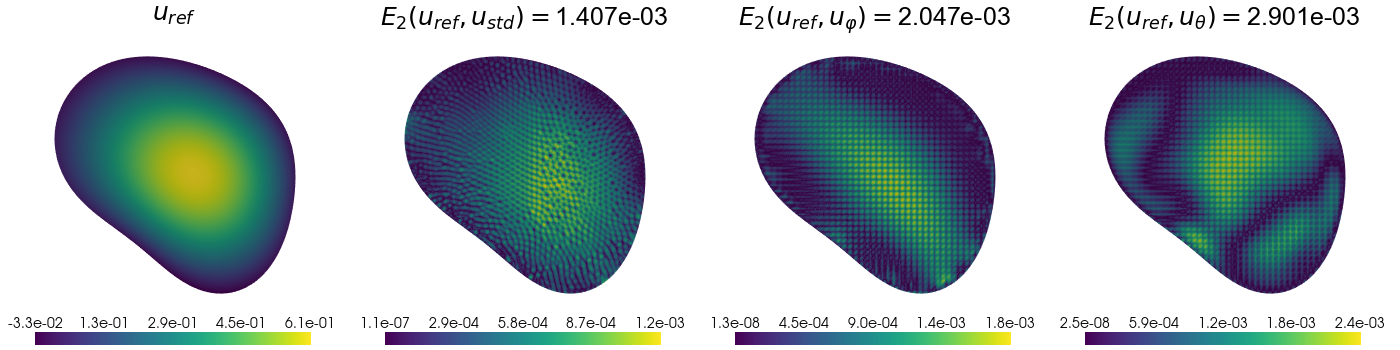}
        \caption{\textbf{Test case 2.} Reference solution ($u_{\text{ref}}$).
            Outputs of standard FEM ($u_{\text{std}}$) , $\varphi$-FEM ($u_{\varphi}$), and ($\varphi$-FEM-FNO $u_{\theta}$).}
        \label{fig:output_random_shapes}
    \end{figure}

    Once again, the results in Fig.\ref{fig:results_random_shapes} (left) demonstrate that $\varphi$-FEM-FNO achieves accuracy comparable to FEM-based methods and outperforms Standard-FEM-FNO. However, $\varphi$-FEM-FNO and Standard-FEM-FNO achieve these results significantly faster, as illustrated in Fig.\ref{fig:results_random_shapes} (right).

    Finally, Fig.~\ref{fig:results_hausdorff_random_shapes} illustrates the correlation between the error and the Hausdorff distance of a test shape to the closest shape in the training data.}

\subsection{2D - Hyperelastic problems with holes}

We now demonstrate the potential of our approach on a test case, close to biomechanical experiments \cite{meunier2008mechanical}:
the non-linear elasticity equation \eqref{eq:governing_hyperelast}. Specifically,
we will consider a rectangular plate with five circular holes, referred to as $\Omega$ and illustrated in Fig.~\ref{fig:plate_with_holes}. Here the variability in the geometry will be the centers {and radii} of the holes.

Let us define the boundaries of the plate as follows:

\begin{itemize}
    \item $\Gamma_D^t$ and $\Gamma_D^b$ are the top and bottom  sides of the rectangle, c.f. Fig.~\ref{fig:plate_with_holes};
    \item $\Gamma_N$ is the Neumann boundary, further subdivided into:
          \begin{itemize}
              \item $\Gamma_N^l$ and $\Gamma_N^r$, are the left and right sides of the rectangle, respectively,
              \item for $i\in\{1,\cdots,5\}$, $\Gamma_N^i$  the boundaries of the $i^{\text{th}}$ circular hole.
          \end{itemize}
\end{itemize}

On $\Gamma_D^b$, the plate is clamped, while a constant displacement {$\mbf u_D$} is applied on $\Gamma_D^t$. For the holes, boundary conditions are enforced using $\varphi$-FEM, whereas standard techniques are used for the remaining boundaries.
\begin{remark}
    We consider the following partition of $\Gamma$:
    \[
        \Gamma = \overbrace{\Gamma_D^b \cup \Gamma_D^t\cup\Gamma_N^l \cup \Gamma_N^r}^{\text{Standard imposition}} \cup  \underbrace{\bigcup_{i=1}^5 \Gamma_N^i}_{\text{$\varphi$-FEM imposition}}\,.
    \]
\end{remark}

The problem is formulated as follows (see \cite[eq. (8.28)]{holzapfel2002nonlinear}): 
find the displacement field $\mbf{u}\in \bb{R}^2$ that satisfies:

\[
    \begin{cases}
        - \Div \mbf{P}(F(\mbf u))              & = 0\,, \quad \hfill \text{in } \Omega \,,           \\
        \hfill \mbf u                          & = \mbf u_D\,, \quad \hfill \text{on } \Gamma_D^t\,, \\
        \hfill \mbf u                          & = 0\,, \qquad \hfill \text{on } \Gamma_D^b\,,       \\
        \hfill \mbf{P}(F(\mbf u)) \cdot \mbf n & = 0\,, \qquad \hfill \text{on } \Gamma_N \,.
    \end{cases}
\]
The first Piola-Kirchhoff stress tensor, $\mbf{P}$, is given by (see \cite[eq. (6.1)]{holzapfel2002nonlinear}):
\[
    \mbf{P}(F(\mbf u)) = \frac{\partial W(F(\mbf u))}{\partial F},
\]
where the strain energy density function $W$ is defined as (see \cite{bonet1997nonlinear}):
\[
    W = \frac{\mu}{2} \left( I_1 - 3 - 2 \ln(J) \right) + \frac{\lambda}{2} \ln(J)^2,
\]
which represents a compressible Neo-Hookean material.

Here, $I_1 = \text{tr}(C)$ is the first invariant of the right Cauchy-Green deformation tensor $C$, defined as $C = F^T  F$, with $F = I + \nabla \mbf{u}$ representing the deformation gradient and $J = \det F$ the Jacobian determinant. The \textit{Lamé} parameters $\mu$ and $\lambda$ are expressed as:
\[
    \mu = \frac{E}{2(1+\nu)}, \text{ and } \lambda = \frac{E \nu}{(1+\nu)(1-2\nu)},
\]
{with the Young modulus of $E$ and the  Poisson's ratio of $\nu$ fixed to $0.97 \, \text{Pa}$ and $0.3$ respectively.}

\subsubsection{$\varphi$-FEM scheme}
In this case, since $\Omega$ is a square domain, the external boundary conditions on $\Omega$ can be applied straightforwardly using conforming methods. Therefore, we will strongly enforce the boundary conditions on all external edges.

To account for the presence of multiple holes, we define distinct level-set functions to construct our $\varphi$-FEM scheme. Each circular hole $\mathcal{C}_i$ with boundary $\Gamma^i_N= \{ \varphi_i =0 \}$, $i=1,\dots, 5$, is defined by
\[ \mathcal{C}_i = \{ \varphi_i < 0 \}, \text{ with } \varphi_i(x,y) =  r_i^2- (x-x_i)^2 - (y-y_i)^2\,, \]
where $(x_i,y_i, r_i)$ are the coordinates of the center and the radius of the hole $i$.

The domain $\Omega$ is then defined as:
\[
    \Omega = \left\{ \underbrace{\prod_{i=1}^5 \varphi_i}_{\varphi} < 0 \right\} \cap (0, 1)^2.
\]

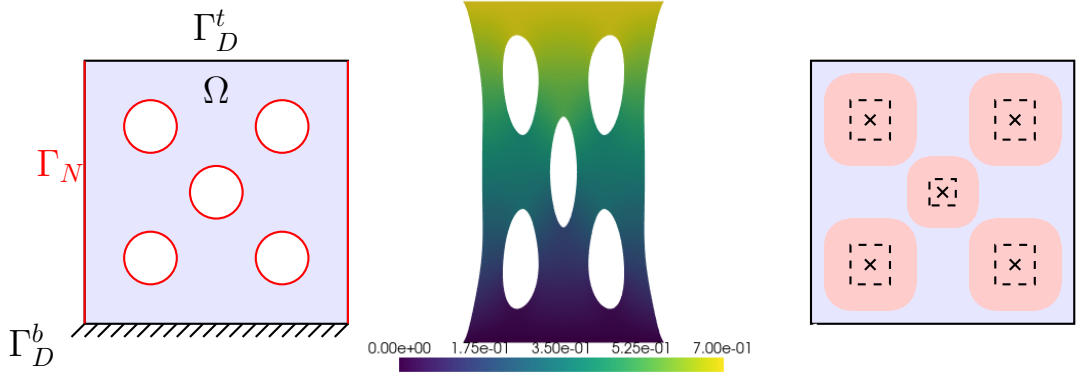
\begin{figure}
    \centering
    \begin{tikzpicture}[scale=3.5]

        \def\rayon{0.1}

        \draw[thick, fill=blue!10] (0,0) rectangle (1,1);
        \draw[thick, red, fill=white] (0.25, 0.75) circle (0.1);
        \draw[thick, red, fill=white] (0.75, 0.75) circle (0.1);
        \draw[thick, red, fill=white] (0.25, 0.25) circle (0.1);
        \draw[thick, red, fill=white] (0.75, 0.25) circle (0.1);
        \draw[thick, red, fill=white] (0.5, 0.5) circle (0.1);

        \draw[thick, red] (0,0) -- (0,1);
        \draw[thick, red] (1,0) -- (1,1);
        \draw(-0.2,-0.2)[above, color = black]node{\Large $\Gamma_D^b$};
        \draw(0.5,1.0)[above, color = black]node{\Large $\Gamma_D^t$};
        \draw(-0.09,0.5)[above, color = red]node{\Large $\Gamma_N$};
        \draw(0.5,0.8)[above, color = black]node{\Large $\Omega$};

        \foreach \i in {0,1,2,...,20} {
                \pgfmathsetmacro\x{\i*0.05}
                \draw[thick, black] (\x,0) -- (\x-0.05,-0.05);
            }

    \end{tikzpicture}
    \includegraphics[width=0.35\textwidth]{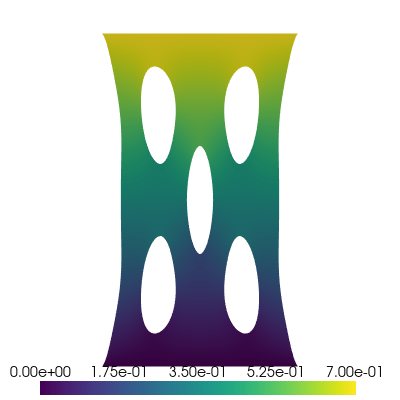}
    \begin{tikzpicture}[scale=3.5]
        \def\rayon{0.1}

        \draw[thick, fill=blue!10] (0,0) rectangle (1,1);

        \fill[fill=red!20] (0.15,0.15) circle (0.1);
        \fill[fill=red!20] (0.3,0.15) circle (0.1);
        \fill[fill=red!20] (0.3,0.3) circle (0.1);
        \fill[fill=red!20] (0.15,0.3) circle (0.1);
        \fill[fill=red!20] (0.05,0.15) rectangle (0.4, 0.3);
        \fill[fill=red!20] (0.15,0.05) rectangle (0.3, 0.4);

        \fill[fill=red!20] (0.7,0.15) circle (0.1);
        \fill[fill=red!20] (0.85,0.15) circle (0.1);
        \fill[fill=red!20] (0.7,0.3) circle (0.1);
        \fill[fill=red!20] (0.85,0.3) circle (0.1);
        \fill[fill=red!20] (0.6,0.7) rectangle (0.95, 0.85);
        \fill[fill=red!20] (0.7,0.6) rectangle (0.85, 0.95);

        \fill[fill=red!20] (0.15,0.7) circle (0.1);
        \fill[fill=red!20] (0.15,0.85) circle (0.1);
        \fill[fill=red!20] (0.3,0.7) circle (0.1);
        \fill[fill=red!20] (0.3,0.85) circle (0.1);
        \fill[fill=red!20] (0.05,0.7) rectangle (0.4, 0.85);
        \fill[fill=red!20] (0.15,0.6) rectangle (0.3, 0.95);

        \fill[fill=red!20] (0.7,0.7) circle (0.1);
        \fill[fill=red!20] (0.7,0.85) circle (0.1);
        \fill[fill=red!20] (0.85,0.7) circle (0.1);
        \fill[fill=red!20] (0.85,0.85) circle (0.1);
        \fill[fill=red!20] (0.6,0.15) rectangle (0.95, 0.3);
        \fill[fill=red!20] (0.7,0.05) rectangle (0.85, 0.4);

        \fill[fill=red!20] (0.45,0.45) circle (0.085);
        \fill[fill=red!20] (0.55,0.45) circle (0.085);
        \fill[fill=red!20] (0.45,0.55) circle (0.085);
        \fill[fill=red!20] (0.55,0.55) circle (0.085);
        \fill[fill=red!20] (0.45,0.45-0.085) rectangle (0.55, 0.55+0.085);
        \fill[fill=red!20] (0.45-0.085,0.45) rectangle (0.55+0.085, 0.55);

        \draw[thick,dashed, fill=red!20] (0.15,0.15) rectangle (0.30, 0.30);
        \draw[thick,dashed, fill=red!20] (0.7,0.15) rectangle (0.85, 0.30);
        \draw[thick,dashed, fill=red!20] (0.15,0.7) rectangle (0.3, 0.85);
        \draw[thick,dashed, fill=red!20] (0.7,0.7) rectangle (0.85, 0.85);
        \draw[thick,dashed, fill=red!20] (0.45,0.45) rectangle (0.55, 0.55);

        \foreach \x/\y in {0.225/0.225, 0.225/0.775, 0.775/0.775, 0.775/0.225, 0.5/0.5} {
                \draw[thick, black] (\x-0.02, \y-0.02) -- (\x+0.02, \y+0.02);
                \draw[thick, black] (\x-0.02, \y+0.02) -- (\x+0.02, \y-0.02);
            }

        \draw(0,-0.2)[above, color = white]node{\Large $\Gamma_D^b$};
    \end{tikzpicture}

    \caption{\textbf{Test case 3.} Left: Representation of the considered situation for the third test case. Center: example of deformed geometry. Right: Representation of the considered variations of the holes for the data generation. The black dashed squares correspond to the bounds of the centers of the holes and the red areas contain all the possible holes configurations. }\label{fig:plate_with_holes}
\end{figure}

An example of this configuration is shown in Fig.~\ref{fig:plate_with_holes}.
{
To account for each boundary condition, we introduce several meshes and sub-meshes. First, we define the computational mesh {$\mathcal{T}_h$}, which covers $\Omega$ and denote by $\Omega_h:=\cup_{T\in \mathcal{T}_h}T$. This mesh consists of all the cells of a Cartesian grid over the box $(0,1)^2$ such that at least one vertex $v$ of the cell satisfies $\varphi(v) < 0$.

Next, we define a sub-mesh $\mathcal{T}_h^{\Gamma}$, which collects all cells intersecting the circular boundaries:

\[ \mathcal{T}_{h}^{\Gamma}:= \{ T \in \mathcal{T}_h : \exists i=1,\dots, 5 \ \text{ s.t. } \ \varphi_i \geqslant 0 \ {\text{ on a vertex of }} T \} \,
\]
and denote by $\Omega_h^{\Gamma}:=\cup_{T\in \mathcal{T}_h^{\Gamma}}T$.
We now introduce the finite element spaces used in the formulation. For an integer $k \geq 2$, the solution $\mbf{u}$ will belong to the space $V_h^k$, defined as:

\[
    V_{h}^{k} := \left\lbrace \mbf{v}_h: \Omega_h \to\mathbb{R}^d : \mbf{v}_{h |T} \in \mathbb{P}^{k}(T)^{d} \ \ \forall T \in \mathcal{T}_h , \ \mbf{v}_h \text{ continuous on } \Omega_h\text{ if } k \ge 0\right\rbrace,
\]
and its homogeneous counterpart $V_{h}^{k,0}$, both FE spaces on $\mathcal{T}_h$.
}
In addition, we need to introduce two auxiliary variables to impose the Neumann boundary conditions on the holes.

Let $\Omega_h^{\Gamma,i}$ be the domain covering the mesh that collects all cells in $\mathcal{T}_h$ cut by the boundary $\Gamma_N^i$:
\[
    \mathcal{T}_h^{\Gamma_i^N} = \{ T \in \mathcal{T}_h : T \cap \Gamma_{i,h}^N \neq \emptyset \},
\]
with $\Gamma_{i,h}^N = \{ \varphi_{i,h} = 0 \}$, where $\varphi_{i,h}$ is the $\mathbb{P}^k$ interpolation of $\varphi_i$ on $\mathcal{T}_h^\Omega$.

The auxiliary variables will live in the following FE spaces:
\[
    Z_{h}^k := \left\lbrace \mbf{z}_h :\Omega_h^\Gamma \to\mathbb{R}^{(d\times d) } : \mbf{z}_{h |T} \in \mathbb{P}^k(T)^{(d\times d)} \ \ \forall T \in \mathcal{T}_h^\Gamma, \ \mbf{z}_h \text{ continuous on }\Omega_h^\Gamma \right\rbrace \,,
\]
and
\[
    Q_{h}^k := \left\lbrace \mbf{q}_h:\Omega_h^\Gamma \to\mathbb{R}^d : \mbf{q}_{h |T} \in \mathbb{P}^{l}(T)^{d} \ \ \forall T \in \mathcal{T}_h^\Gamma , \ \mbf{q}_h \text{ continuous on }\Omega_h^\Gamma\text{ if }k\ge 0\right\rbrace.
\]

For each hole $i$, we impose homogeneous Neumann boundary conditions through the following equations:

\begin{align*}
    \mbf{y} + \mbf{P}(F(\mbf{u}))=0,            & \quad \text{on } \Omega_h^{\Gamma,i}\,, \\
    \mbf{y}\nabla\varphi_i+ \mbf{p}\varphi_i=0, & \quad \text{on } \Omega_h^{\Gamma,i}\,.
\end{align*}

The variables $\mbf{y}$ and $\mbf{p}$ are discretized in the spaces $Z_h^k$ and $Q_h^{k-1}$, respectively.

This yields to the following variational formulation: find $\mbf{u}_{h} \in V_{h}^k$, $\mbf{p}_h\in Q_{h}^{k-1}$, $\mbf{y}_{h}\in Z_{h}^k$, such that

\begin{multline}\label{eq:scheme_plate}
    \int_{\Omega_{h}} \mbf{P} (F(\mbf{u}_h)) : \nabla \mbf{v}_h
    + \sum_{i=1}^{5} \bigg( \int_{\partial\Omega_h^{\Gamma, i}} \mbf{y}_h \mbf{n} \cdot \mbf{v}_h + \gamma_{u} \int_{\Omega_h^{\Gamma, i}} (\mbf{y}_{h} + \mbf{P} (F(\mbf{u}_{h}))) : (\mbf{z}_{h} + D_{\mbf{u}}(\mbf{P}\circ F )(\mbf{u}_{h})\mbf{v}_{h}) \\
    + \frac{\gamma_{p}}{h^2} \int_{\Omega_h^{\Gamma,i}} (\mbf{y}_{h} \nabla \varphi_{i,h} + \frac{1}{h} \mbf{p}_{h} \varphi_{i,h}) \cdot (\mbf{z}_{h} \nabla \varphi_{i,h} + \frac{1}{h} \mbf{q}_{h} \varphi_{i,h}) \\
    + \gamma_{div} \int_{\Omega_h^{\Gamma, j}}\Div \mbf{y}_h \cdot \Div \mbf{z}_h \bigg)
    + G_h \left( \mbf{u}_{h}, \mbf{v}_{h} \right)
    =  0\,, \\
    \forall \mbf{v}_{h} \in V_{h}^{k,0}, \ \mbf{q}_h\in Q_{h}^{k-1}, \ \mbf{z}_{h}\in Z_{h}^k\,,
\end{multline}
where
\begin{equation*}
    G_h(\mbf{u}, \mbf{v}) : = \sigma_N h \int_{\Gamma_h} \left[ \mbf{P}(F(\mbf{u}))\mbf{n}  \right] \cdot \left[ D_{\mbf{u}}(\mbf{P}\circ F)(\mbf{u})\mbf{v}\mbf{n}\right]\,,
\end{equation*}
with $\Gamma_h:=\partial \Omega_h^{\Gamma}\setminus\partial\Omega_h$ (corresponding to the facets between the blue and grey cells in the example of Fig.~\ref{fig:domain_submeshes}), $D_{\mbf{u}}(\mbf{P}\circ F)(\mbf{u})\mbf{v}$ denoting the derivative of $\mbf{P}$ evaluated at $\mbf{u }$, in the direction $\mbf v$ and $\gamma_p$, $\gamma_u$, $\gamma_{div}$, $\sigma_N$ some positive constants.

\begin{remark}
    The third term in \eqref{eq:scheme_plate} is the differential at $(\mbf{u}_{h},\mbf{y}_{h})$ in the direction $(\mbf{v}_{h},\mbf{z}_{h})$ of
    $$    \int_{\Omega_h^{\Gamma, i}} (\mbf{y}_{h} + \mbf{P}(F (\mbf{u}_{h}))):(\mbf{y}_{h} + \mbf{P}(F(\mbf{u}_{h}))).$$
    Hence,  it is a quid of penalization associated with the constraint
    $$\mbf{y}_{h} =- \mbf{P}(F (\mbf{u}_{h})).$$
\end{remark}

\begin{figure}
    \centering
    \includegraphics[width=0.45\textwidth]{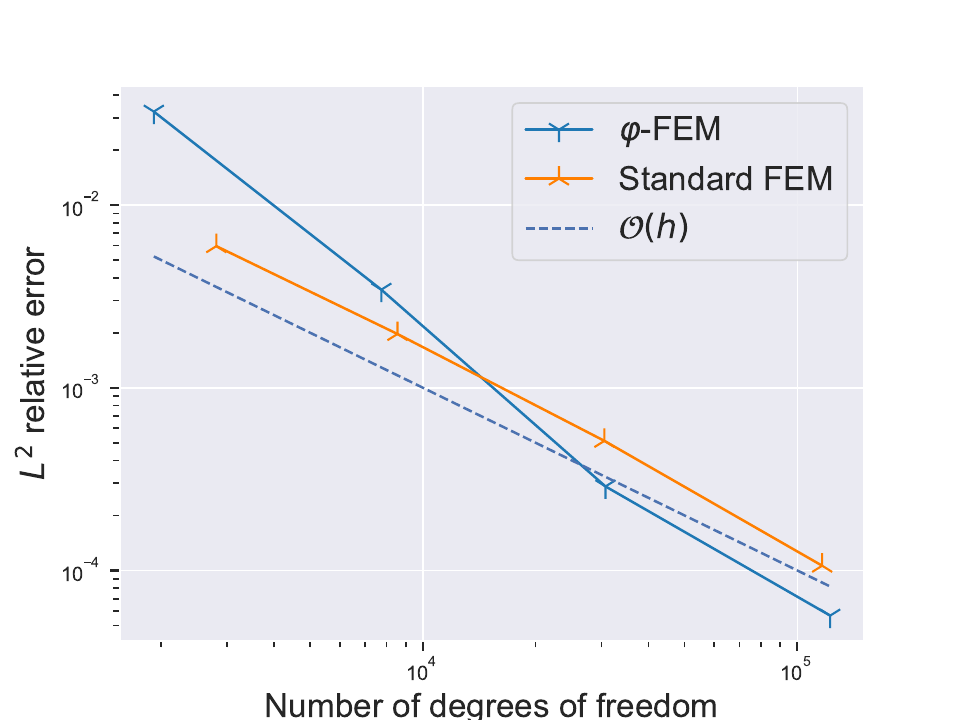}
    \includegraphics[width=0.45\textwidth]{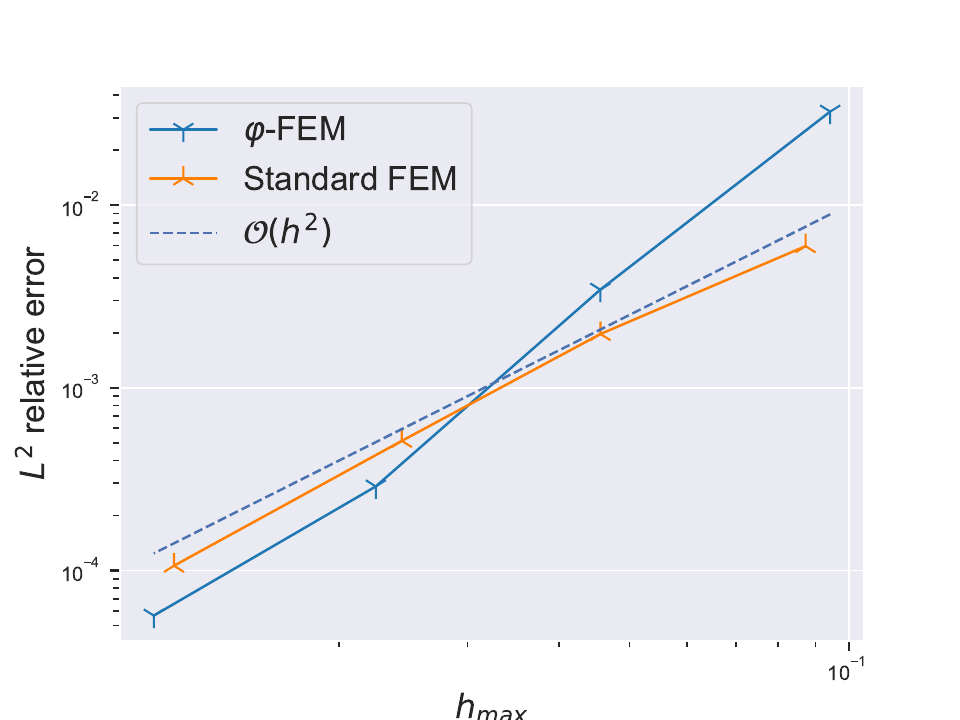}
    \caption{\textbf{Test case 3.} Convergence of the two finite element methods, with respect to the number of degrees of freedom (left) and the maximal size of cell (right).}
    \label{fig:convergence_fems_plate}
\end{figure}

We validate this $\varphi$-FEM scheme by comparing the convergence of the method with the convergence of a standard finite element method, computing the relative $L^2$ error with a reference fine standard FEM solution. Referring to Fig.~\ref{fig:convergence_fems_plate}, we can conclude that $\varphi$-FEM outperforms the standard method for sufficiently fine meshes.

\subsubsection{$\varphi$-FEM-FNO results}
{Since this test case is more complex, particularly due to the treatment of Neumann boundary conditions using the $\varphi$-FEM approach, the $\varphi$-FEM-FNO architecture presented earlier needs slight modifications. Specifically, in this approach, the solution is not multiplied by the level-set function over the entire domain but only on the boundary cells. For more details on this method applied to the Poisson equation, refer to \cite{phiFEM2}, or for the linear elasticity equation, see \cite{cotin:hal-03372733}. As a result, the neural operator no longer predicts
$w$; instead, it directly provides the
$\varphi$-FEM solution to the problem. This modified method is referred to as $\varphi$-FEM-FNO-2 in the first test case.

Moreover, as this problem involves only null right-hand sides except for $\mbf{u}_D$, the ground-truth operator to approximate is defined as:

\begin{equation*}
    \begin{array}{rccl}
        \mathcal{G}^{\dagger}: & \mathbb R^{n_x\times n_y \times 2} & \to     & \mathbb R^{n_x\times n_y\times 2} \\
                               & (\varphi_h,g_{h,y})                & \mapsto & \mbf{u}_h = (u_{h,x}, u_{h,y})\,,
    \end{array}
\end{equation*}
where $u_{h,x}$ and $u_{h,y}$ are the two components of the solution vector $\mbf{u}_h$, and $g_{h,y}$ represents the vertical component of the Dirichlet boundary condition $\mbf{u}_D$, which is constant throughout the domain, i.e., $g_{h,y} = g$ for every pixel.

For this test case, the neural network approach demonstrates significant advantages over FEM-based methods. The problem's strong non-linearity typically requires iterative solvers in classical methods, and incremental forces are often applied to prevent divergence. These iterative processes are computationally expensive. In contrast, the $\varphi$-FEM-FNO approach only requires a well-trained operator to directly obtain the solution. Although data generation is more time-consuming for this test case compared to simpler ones, the trained operator enables extremely fast solution computation.

\paragraph{Data generation}
To generate training, validation, and test data, we consider a configuration with five circular holes. The holes are placed sufficiently far from the boundaries of the unit box and positioned to avoid interpenetration. Data generation is performed using a Latin Hypercube sampling strategy, this time in a 16-dimensional space (15 dimensions for the hole parameters and one for the applied boundary condition at the top of the box). A graphical representation of the sampled hole configurations is shown in Fig.~\ref{fig:plate_with_holes} (right).
The parameters of the $\varphi$-FEM scheme are set to $\gamma_u = 0.001$, $\gamma_p = \gamma_{div} = \sigma_N = 0.01$.

\paragraph{Loss modification}

To train the operator, we need to adapt the loss function defined in \eqref{eq:loss}. We will now use only an approximation of the $H^1$ semi-norm as loss function, defined
by
\[
    \mathcal{L}\left(U_{\text{true}}; U_{\theta}\right) = \frac{1}{N_\text{data}}\sum_{n=0}^{N_\text{data}}
    \left(\mathcal{E}_1( u^n_{\text{true}, x} ; u^n_{\theta, x}) + \mathcal{E}_1( u^n_{\text{true}, y} ; u^n_{\theta, y}) \right)\,,
\]
where
\[
    \mathcal{E}_1(u^n_{\text{true}, \cdot}; u^n_{\theta,  \cdot}) =
    \| \nabla_x^h u^n_{\text{true},  \cdot} - \nabla_x^h u^n_{\theta,  \cdot} \|^2_{0,\mathcal{S}_1^n}
    + \| \nabla_y^h u^n_{\text{true},  \cdot} - \nabla_y^h u^n_{\theta,  \cdot} \|^2_{0, \mathcal{S}_1^n} \,.
\]

\begin{remark}
    Using the $H^1$ semi-norm instead of the full $H^1$ norm enhances the performance of the operator, particularly in the application of boundary conditions. Once the operator generates a prediction, we can "adjust" the solution by subtracting the mean value of the prediction at the lower boundary nodes, where the solution is known to be zero.

    This approach offers several advantages. First, it simplifies the optimization process, as the loss function involves fewer terms. Second, it reduces the error at the boundary nodes compared to using the full $H^1$ norm, leading to more accurate boundary condition enforcement.
\end{remark}
\paragraph{Results}

To measure the performance of the finite element based approaches and of the $\varphi$-FEM-FNO approach, we will compute the $L^2$ relative error, between a reference displacement $\mbf{u}_\text{ref}$ and an approximation $\mbf{u}_h$, denoted $\bar{L_2}(\mbf{u}_\text{ref}, \mbf{u}_h)$.

We train the operator for 2000 epochs using 200 training data divided into batches of size 8 and 300 validation data. We then evaluate the performance of the trained operator on 300 test data, compared to a standard finite element method and to $\varphi$-FEM both using triangular $\bb{P}^2$ elements, with mesh sizes such that the total number of degrees of freedom is close to the dimension of the images used for $\varphi$-FEM-FNO (i.e. $2 \times n_x \times n_y$).
We represent in Fig.~\ref{fig:outputs_plaque} an example of displaced mesh, and the difference between the approximated results of the methods and a reference fine FEM solution.

\begin{figure}
    \centering
    \includegraphics[width=\textwidth]{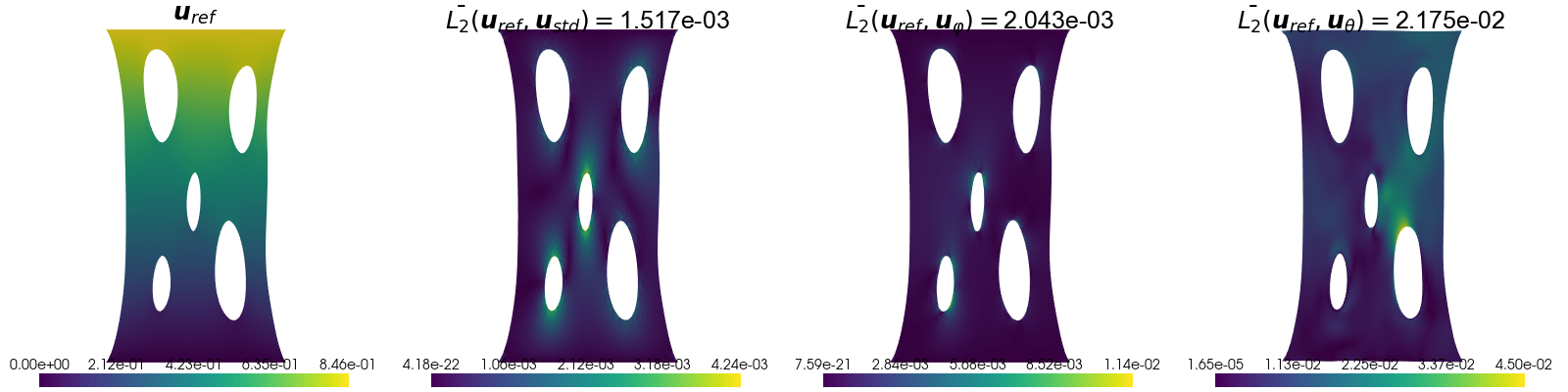}
    \caption{\textbf{Test case 3}. Example of output of the different methods, corresponding to the median of the relative $L^2$ errors of $\varphi$-FEM-FNO among the 300 problems.}
    \label{fig:outputs_plaque}
\end{figure}

In addition, in Fig.~\ref{fig:results_plaque} (left), we compare the $L^2$ errors of the 3 methods on the 300 test problems. In Fig.~\ref{fig:results_plaque} (right), we give the relative Hausdorff error for the same 300 problems, with respect to the computation time. These representations illustrate that the FEM's based approaches are only at most 10 times more precise than our approach, while the computation times are close to 1000 times higher for the FEM's based approaches.

\begin{figure}
    \centering
    \includegraphics[width=0.49\textwidth]{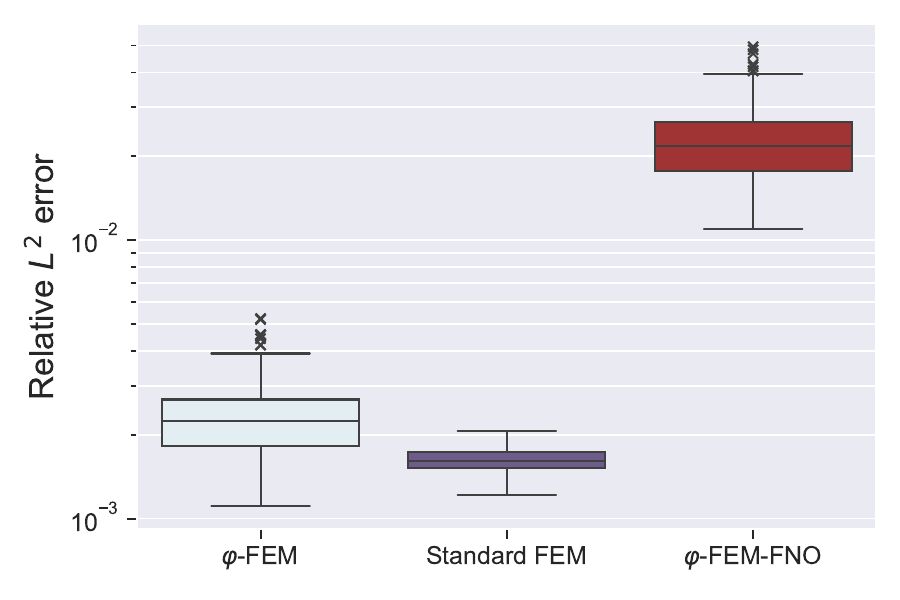}
    \includegraphics[width=0.49\textwidth]{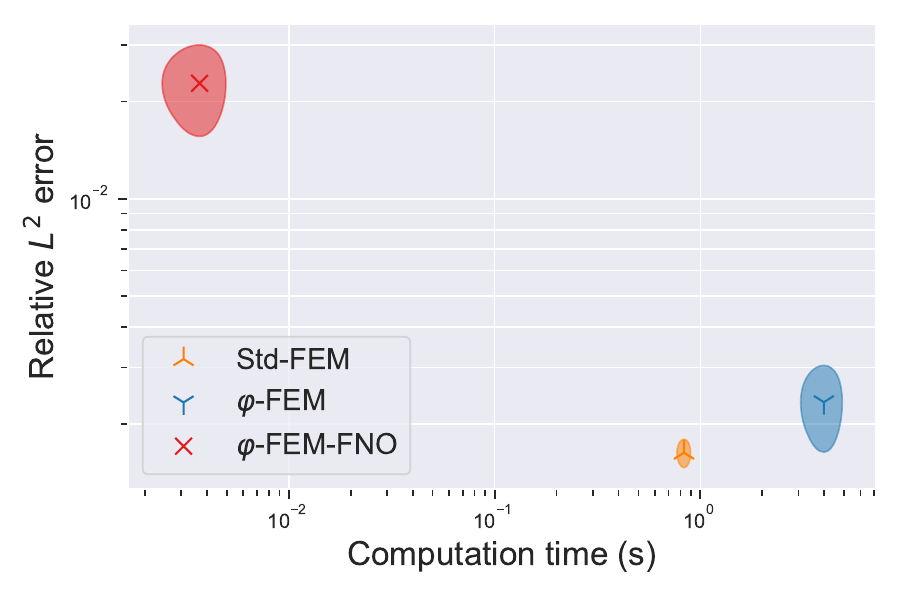}
    \caption{\textbf{Test case 3}. Left: $L^2$ errors of the methods. Right: relative $L^2$ errors of the methods against computation time.}
    \label{fig:results_plaque}
\end{figure}
}

\section{Conclusion and future works}\label{sec:concl}

We have shown on three test cases that after training, our $\varphi$-FEM-FNO can compute faster than standard finite element methods, $\varphi$-FEM, an interpolate-FEM approach, or Geo-FNO on several problems. Moreover, we have illustrated that these results can be obtained using small amounts of training data, even for complex cases with big variations of geometries or non-linear equations.

A number of perspectives remain for future research. It would be interesting to extend the results to other problems since $\varphi$-FEM schemes have been written and studied theoretically and numerically (mixed conditions, Stokes, time-dependent PDEs, \dots).

Moreover, we have introduced and validated a new $\varphi$-FEM scheme to treat a case of non-linear elastic equation. However, in the future, we can extend our results to the case of other hyperelastic materials as in \cite{odot:hal-03327818}, and implement the method in the DeepPhysX project\footnote{\url{https://mimesis.inria.fr/project/deepphysx/}}.
Furthermore, another interesting point would be to extend our method to more realistic scenarios, considering real medical images and more realistic forces and boundary conditions. Finally, we can also imagine extending our method to the case of $\bb{P}^k$ functions, using the degrees of freedom values instead of nodal values for the data generation and thus predicting the values of the solution at each $\bb{P}^k$ degrees of freedom.

Finally, to represent more complex and general forces, one can train an FNO using Gaussian forces. Then, one can decompose a new random force in a sum of Gaussian distributions and use the trained model on each one of the Gaussian forces. Thanks to GPU parallelization, each prediction can be done simultaneously, and it only remains to sum the predictions to obtain the final result.

\section{ Acknowledgment }

The authors were supported by the ANR project JCJC 22-CE46-0003. The authors would like to thank Nicola Zotto and Sidaty El Hadramy for their remarks and help during the preparation of this paper.

\bibliographystyle{abbrv}
\bibliography{biblio}

\appendix

\section{FNO implementation details}
\subsection{Standardization of the data}\label{sec:norm}
To improve the performance of our FNO, since the data can have very different values, we have decided to standardize the input and output data, as in \cite{paper_FNO}.
The standardization is applied independently channel by channel of $X$.
For each channel $C$ of $X$, denoting by $C^{\text{train}}$ the training part of the data-set corresponding to the channel $C$, the associated standardized channel is given by
\[
    N_C(C) = \left( \frac{C - \text{mean} (C^{\text{train}})}{\text{std} (C^{\text{train}})} \right)\,,
\]
where the mean and standard-deviation are computed only on $\Omega_h$, since all the values are 0 outside $\Omega_h$.

The unstandardization function $N^{-1}$ is given by
\[
    N^{-1}(Y) = Y \times \text{std} (Y^{\text{train}})+ \text{mean}(Y^{\text{train}})\,,
\]
where $Y$ denotes a channel of the output of the FNO and $Y^{\text{train}}$ is the vector composed of the training ground truth solutions.

\subsection{ADAM and training loop algorithm}\label{Ap:algo}
We consider here the case of the first test case. We present the details of the considered ADAM optimizer in Algorithm~\ref{algo:ADAM}.
In Algorithm~\ref{algo:training} we denote $(F^i, \varphi^i, G^i)$ a batch of data. The batches are randomly chosen such that $F^i = (f_h^k)_{k \in K_i}$, $\varphi^i = (\varphi_h^k)_{k \in K_i}$, $G^i = (g_h^k)_{k \in K_i}$ with $K_i$ a collection of random indices of data and $i \in \{1, \dots, \text{number of batches}\}$. The sets $K_i$ are constructed such that $K_i \cap K_j = \emptyset$ for $i\neq j$.
\begin{algorithm}
    \caption{ADAM optimizer step.}\label{algo:ADAM}
    \begin{algorithmic}
        \STATE Initialisation : $t$, $\theta_{t-1}$, $\beta_1$, $\beta_2$, $\varepsilon$, $m_{t-1}$, $v_{t-1}$.
        \STATE Compute the gradient : $g_t \gets \nabla f(\theta_{t-1})$
        \STATE Momentum update :
        \[ m_t \gets \beta_1 \cdot m_{t-1} + (1 - \beta_1) \cdot g_t\,, \qquad
            v_t \gets \beta_2 \cdot v_{t-1} + (1 - \beta_2) \cdot g_t \cdot \bar{g_t} \]
        \vspace{-0.4cm}
        \STATE Bias correction : \[ \hat{m}_t \gets \frac{m_t}{1 - \beta_1^t}\,,
            \qquad \hat{v}_t \gets \frac{v_t}{1 - \beta_2^t}\] \vspace{-0.4cm}
        \STATE Parameters update :
        \[ \theta_t \gets \theta_{t-1} - \frac{\alpha}{\sqrt{\hat{v}_t} + \varepsilon} \cdot \hat{m}_t - w_1 \theta_{t-1} \]
    \end{algorithmic}
\end{algorithm}

\begin{remark}[Calibration of the learning rate.]
    The learning rate is a critical parameter to tune for achieving accurate results. While we have not included specific results to illustrate our choice of learning rate, extensive testing was conducted to determine the optimal value. Choosing a learning rate that is too high or decreasing it too slowly results in significant oscillations and poor convergence. Conversely, selecting a learning rate that is too low or decreasing it too quickly leads to slow and suboptimal convergence, as the loss decreases very slowly and fails to reach sufficiently low values to produce good results.

    To address this, we fine-tuned the learning rates through multiple training sessions on both test cases, experimenting with various learning rate schedulers. The scheduler that provided the best results was selected, using the validation loss as a criterion to adjust the learning rate dynamically.

\end{remark}
\begin{algorithm}
    \caption{Training loop.}\label{algo:training}
    \begin{algorithmic}
        \STATE \textbf{Initialisation:} $\theta_0$ the initial random parameters, $X = (F, \varphi, G)$ and $Y_{\text{true}}$ the training part of the dataset, the batch size and the regularization parameter $\lambda$ .
        \FOR{$t=1$ \TO number of epochs}
        \FOR{$i=1$ \TO number of batches \tikzmark{top_training}}
        \STATE Select a batch $ (F^i, \varphi^i, G^i) \subset X$ and $Y_\text{true}^i \subset Y_\text{true}$ of size batch size.
        \STATE Call the model :
        $Y_{\theta} =\mathcal{G}_{\theta_{ti-1}}(F^i, \varphi^i, G^i)$.
        \STATE Compute the loss : \[\mathcal{L}(Y_{\text{true}}^i, Y_{\theta}) + \underbrace{\frac{\lambda}{2 \times \text{batch size}}\sum_j | w_j |^2}_{\text{$L^2$ regularization}}\,. \]
        \STATE Compute the gradient of the loss, w.r.t the parameters $\theta_{ti-1}$: $\nabla_{\theta_{ti-1}} \mathcal{L}$.
        \STATE Optimizer step : step of Algorithm \ref{algo:ADAM}.
        \ENDFOR \tikzmark{bottom_train}
        \vspace{0.2cm}
        \STATE Let $(F_\text{val}, \varphi_\text{val}, G_\text{val})$ and $Y_{\text{val}}$ be the validation part of the dataset. \tikzmark{top_validate}\hspace{2cm} \tikzmark{right}
        \STATE Call the model on the validation sample : $Y_{\theta} = \mathcal{G}_{\theta_{ti}}(F_\text{val}, \varphi_\text{val}, G_\text{val})$.
        \STATE Compute the loss : $\mathcal{L}(Y_{\text{val}}, Y_{\theta})$. \tikzmark{bottom_validate}
        \STATE Learning rate scheduler step.
        \ENDFOR
    \end{algorithmic}
    \AddNote{top_training}{bottom_train}{right}{black}{\qquad Training step}
    \AddNote{top_validate}{bottom_validate}{right}{black}{\qquad Validation step}
\end{algorithm}

\end{document}